\documentclass[12pt,a4paper]{article}
\usepackage{amssymb,amsmath,amsthm}
\usepackage{latexsym}
\usepackage{color}
\usepackage{stmaryrd}

\newtheorem{theorem}{Theorem}
\newtheorem{lemma}{Lemma}
\newtheorem{proposition}{Proposition}
\newtheorem{corollary}{Corollary}
\newtheorem{remark}{Remark}
\numberwithin{equation}{section}
\numberwithin{theorem}{section}
\numberwithin{lemma}{section}
\numberwithin{proposition}{section}
\numberwithin{corollary}{section}
\numberwithin{remark}{section}

\begin{document}
\title{Weighted $L^q$-$L^r$ estimates of a generalized Oseen evolution operator in exterior domains, with applications to the Navier-Stokes flow past a rigid body}
\author{Yosuke Asami\\
Graduate School of Mathematics \\
Nagoya University \\
Nagoya 464-8602, Japan \\
\texttt{yosuke.asami.e6@math.nagoya-u.ac.jp} \\
}
\date{}
\maketitle
\begin{abstract}
We consider a generalized Oseen evolution operator arising from time-dependent rigid motions in 3D exterior domains.
Hishida and the present author established the regularity theory of the evolution operator by use of the $L^q$-$L^r$ smoothing estimates with a weight function $(1+|x|)^{\alpha}$.
The objective of this paper is to develop the weighted $L^q$-$L^r$ estimates globally in time. 
These estimates are applied to the study of the large time behavior of the Navier-Stokes flow around a time-periodic solution in weighted spaces.


\end{abstract}

\section{Introduction}
We consider the motion of a viscous incompressible fluid filling the exterior of an obstacle in $\mathbb R^3$, where the obstacle is a rigid body whose motion is prescribed and can be time-dependent.
By taking a frame attached to the moving body, we can reduce the problem in the inertial frame to an equivalent one in a time-independent exterior domain $D$, whose boundary $\partial D$ is of class $C^{1, 1}$.
Then the unknown fluid velocity field $u=u(x, t)\in \mathbb R^3$ and the scalar pressure $p_u=p_u(x, t)\in \mathbb R$ obey the following Navier-Stokes system (see \cite{Ga02}):
\begin{equation}
\begin{split}
&\left.
\begin{array}{rl}
\partial_tu+u\cdot\nabla u&=\Delta u+(\eta+\omega\times x)\cdot\nabla u-\omega\times u-\nabla p_u \\
\mbox{div $u$}&=0 \\
\end{array}
\right\} \;\;\mbox{in $D\times (0,\infty)$},  \\
&\quad u|_{\partial D}=\eta+\omega\times x, \qquad
\lim_{|x|\to\infty}u=0, \qquad
u(\cdot,0)=u_0.
\end{split}
\label{ns}
\end{equation}
The functions $\eta=\eta(t)\in \mathbb R^3$ and $\omega=\omega(t)\in\mathbb R^3$ denote given translational and angular velocities of the rigid body, respectively, while $u_0=u_0(x)\in \mathbb R^3$ stands for a given initial velocity.
Here and hereafter, $\times$ denotes the exterior product.
In this paper, we deduce weighted estimates of the solution to the linearlized system for \eqref{ns} globally in time and, as an application, study the large time behavior of the Navier-Stokes flow around a specific time-periodic motion whose summability at spatial infinity is scale-subcritical.
\par
The well-posedness of the initial value problem \eqref{ns} is challenging because of the non-autonomous character as well as a hyperblic aspect caused by the drift term $(\omega\times x)\cdot\nabla u$.
Regarding the latter, indeed, even for the autonomous case, the corresponding semigroup generated by the linearized system \eqref{linear} below is not an analytic one (\cite{FN07, Hi99}) unless $\omega\equiv 0$.
Nevertheless, there are notable existence theorems such as those by \cite{GS05, GHH06, Hi99, HiShi}, see \cite{Hi13} and the references therein.
In what follows, we focus on the problem within the non-autonomous regime, that is, the rigid motion $\eta+\omega\times x$ is actually time-dependent. 
Hansel and Rhandi \cite{HR14} proved the existence of a unique mild solution to \eqref{ns} locally in time under the assumption $\eta, \omega\in C^1([0, \infty); \mathbb R^3)$ when the initial velocity satisfies
\begin{align*}
u_0\in L^q(D)\;\;\mbox{with}\;\; q\in [3, \infty),\qquad \nu\cdot(u_0-\eta(0)-\omega(0)\times x)|_{\partial D}=0,\qquad \mathrm{div}\,u_0=0,
\end{align*}
where $\nu$ stands for the outer unit normal vector to $\partial D$.
What is remarkable is that the authors of \cite{HR14} constructed the evolution operator $T(t, s)$ in $L^q_{\sigma}(D)$ (solenoidal $L^q$-space) for every  $q\in (1, \infty)$, which is a solution operator $u(s)=f\mapsto u(t)=T(t, s)f$ to the linearized system 
\begin{equation}
\begin{split}
&\partial_tu+{\mathcal L}(t)u
+\nabla p_u=0, \qquad \mbox{div $u$}=0 \quad\mbox{in $D\times (s,\infty)$}, \\
&u|_{\partial D}=0, \qquad \lim_{|x|\to\infty}u=0, \qquad  u(\cdot,s)=f,
\end{split}
\label{linear}
\end{equation}
under the assumption
\begin{align}\label{ass1-rigid}
\eta, \omega\in C^{0, \theta}_{\mathrm{loc}}([0, \infty); \mathbb R^3)
\end{align}
for some $\theta\in (0, 1]$ without relying on any general theory, and deduced the $L^q$-$L^r$ smoothing estimates near the initial time,  where $\mathcal L(t)$ is the differential operator given by
\begin{align}
{\mathcal L}(t)=-\Delta -(\eta(t)+\omega(t)\times x)\cdot\nabla+\omega(t)\times.
\label{pre-proj-0}
\end{align}
On account of the hyperbolic aspect mentioned above, however, the evolution operator $T(t, s)$ is not of parabolic type in the sense of Tanabe-Sobolevski\u\i,\,see \cite{Ta}.
Thus, the regularity properties of $T(t, s)$ such as estimate of temporal derivative as well as the H\"older estimate are not clear in \cite{HR14}, and it is not known whether the solution obtained in \cite{HR14} is a strong one.
Recently, Hishida and the present author \cite{AH25} showed the existence of a strong $L^q$-solution to \eqref{ns} in the sense that it is of class $C^1$ in time with values in $L^q(D)$ when assuming 
\begin{align*}
\rho u_0 \in L^q(D)\;\;\mbox{with}\;\; q\in (3, \infty),\qquad \nu\cdot(u_0-\eta(0)-\omega(0)\times x)|_{\partial D}=0,\qquad \mathrm{div}\,u_0=0,
\end{align*}
as well as appropriate regularity of $\eta$ and $\omega$, where $\rho$ is the weight function given by
\begin{align}\label{weight}
\rho=\rho(x):=1+|x|,
\end{align}
which is related to the drift term $(\omega\times x)\cdot\nabla u$.
We emphasize that this result is new even for the autonomous case.
One of key ingredients is to show the weighted $L^q$-$L^r$ smoothing estimates of $T(t, s)f$ with the initial data $f\in L^q_{\sigma}(D)$ satisfying 
\begin{align}\label{intro-wei}
\rho^{\alpha} f\in L^q(D)\quad\mbox{for}\quad 0\le\alpha<3,\quad \frac{3}{3-\alpha}<q<\infty,
\end{align}
see the item 1 of Proposition \ref{prop1-gen-sm} below.
With the aid of these estimates, it was proved in \cite{AH25} that the evolution operator $T(t, s)f$ provides a strong solution to \eqref{linear} for the initial data $f\in L^q_{\sigma}(D)$ fulfilling $\eqref{intro-wei}$ with $\alpha=1$.
As for the global mild solution to \eqref{ns}, it is possible to discuss various situations \cite{GH21, GH26, Huy23, T22, T26} because the $L^q$-$L^r$ decay estimates of the evolution operator $T(t, s)$ were established well by Hishida \cite{Hi18, Hi20} under the assumption
\begin{align}\label{ass2-rigid}
\eta, \omega\in C^{0, \theta}([0, \infty); \mathbb R^3)\cap L^{\infty}(0, \infty; \mathbb  R^3)
\end{align}
for some $\theta\in (0, 1]$, see also Takahashi and Watanabe \cite{TW} even for the case of the Lipschitz boundary.
In view of the aforementioned theory developed by \cite{AH25}, it seems crucial to deduce the weighted $L^q$-$L^r$ estimates of $T(t, s)f$ globally in time with $f\in L^q_{\sigma}(D)$ satisfying \eqref{intro-wei} for the case $\alpha=1$ in order to obtain a global strong solution to \eqref{ns}.
\par
Due to the inhomogeneous boundary condition $u|_{\partial D}=\eta+\omega\times x$, we should take some nontrivial solution as a basic flow, around which we are going to construct a global solution to \eqref{ns}.
Since the system is non-autonomous, it is reasonable to consider a time-periodic flow $(v, p_v)$, which satisfies the problem
\begin{align}\label{ns-per-int}
\begin{split}
&\left.
\begin{array}{rl}
\partial_tv+v\cdot\nabla v&=\Delta v+(\eta+\omega\times x)\cdot\nabla v-\omega\times v-\nabla p_v \\
\mbox{div $v$}&=0 \\
\end{array}
\right\} \;\;\mbox{in $D\times (-\infty,\infty)$},  \\
&\quad v|_{\partial D}=\eta+\omega\times x, \qquad
\lim_{|x|\to\infty}v=0,
\end{split}
\end{align}
where $\eta$ and $\omega$ are assumed to be $\mathcal T$-periodic with some $\mathcal T>0$.
Galdi \cite{Ga22} and Tomoki Takahashi \cite{T26} proved the unique existence of a periodic solution if $\mathcal T$-periodic rigid motion $\eta+\omega\times x$ satisfies either of the assumptions
\begin{align}
&\mbox{(i) $\eta(t)=0$ or $\omega(t)=0$;}\label{ass1-per}\\
&\mbox{(ii) $\eta(t)$ and $\omega(t)$ are parallel to the same constant vector,}\label{ass2-per}
\end{align}
as well as the smallness of $\eta, \omega\in W^{2, 2}(0, \mathcal T; \mathbb R^3)$, see Proposition \ref{prop-galtak}.
Their solutions also enjoy the spatial pointwise behavior reflecting the wake structure.
To describe their results, we set 
\begin{align}\label{def-lambda}
\lambda e_1:=\frac{1}{\mathcal T}\int_{0}^{\mathcal T}\eta(\tau)\,d\tau,
\end{align}
where the average of the translational velocity over the period is assumed to be parallel to $e_1$ with $\lambda\ge 0$ without loss of generality.
If in particular $\lambda>0$, then the time-periodic flow $v$ satisfies the pointwise estimates
\begin{align}\label{subcrit}
\sup_{t\in\mathbb R}|\nabla^jv(x, t)|\le \frac{C}{(1+|x|)^{1+\frac{j}{2}}(1+\lambda|x|+\lambda x_1)^{1+\frac{j}{2}}}\|(\eta, \omega)\|_{W^{2, 2}(0, \mathcal T)},
\end{align}
for $j=0$ (\cite{Ga22}) and $j=1$ (\cite{T26}).
Note that the estimate \eqref{subcrit} implies that the periodic flow $v(t)$ possesses so fine summability at infinity that it is of scale-subcritical class: $v(t)\in L^{q_0}(D)$ with some $q_0<3$ and $\nabla v(t)\in L^{q_1}(D)$ with some $q_1<3/2$, see Lemma \ref{lem-stab-sum}.
Successively, in \cite{T26}, among others, stability of this periodic flow was proved within the $L^q$-framework.
Let us also mention that the pointwise estimate \eqref{subcrit} with $j=0$ (in which $W^{2, 2}(0, \mathcal T; \mathbb R^3)$ is replaced by $W^{1, \infty}(0, \mathcal T; \mathbb R^3)$) was derived by Hishida \cite{Hi22} under rather comprehensive assumptions, which cover \eqref{ass1-per}--\eqref{ass2-per}, through a different approach from theirs, see \cite[Remark 2.3]{Hi22}. 
Tomoki Takahashi \cite[Remark 3.8]{T26} also deduced the pointwise behavior \eqref{subcrit} with $j=1$ under the same assumptions as in \cite{Hi22}.
\par
In this paper, as the basic flow, we take the time-periodic solution fulfilling \eqref{subcrit} with $\lambda>0$ obtained in \cite{Ga22, T26}.
Our aim is to develop the weighted $L^q$-$L^r$ estimates of the evolution operator $T(t, s)$ globally in time and then to derive the large time behavior of a global solution to \eqref{ns} around the periodic solution above in weighted spaces by applying those estimates of the evolution operator.
To this end, we follow the argument adopted by Tomoki Takahashi \cite{T24}, who studied a global solution to \eqref{ns} with anisotropic weights $(1+|x|)^{\alpha_1}(1+|x|-x_1)^{\alpha_2}$ reflecting the wake structure when the constant vector $\eta$ is parallel to $e_1$ and $\omega\equiv 0$.
The essential step of his analysis is to deduce the weighted $L^q$-$L^r$ estimates of the Oseen semigroup globally in time with such weights, where those estimates are inhomogeneous in the sense that the weighted $L^r$-norm of the semigroup is bounded from above by $Ct^{-\kappa}\|\rho^{\alpha}f\|_{q}+Ct^{-\kappa}(1+t)^{\alpha}\|f\|_q$ if, in particular, the weight is given by $\rho^{\alpha}$, where $\kappa=(3/q-3/r)/2$ and $\|\cdot\|_{q}$ denotes the norm of the space $L^q(D)$.
An important observation of \cite{T24} is that the linear estimate of this form is useful in order that estimates are closed for the nonlinear problem.
\par
Our first theorem (Theorem \ref{thm1} below) is the weighted $L^q$-$L^r$ estimates of the evolution operator $T(t, s)$ globally in time under the assumption \eqref{ass2-rigid} for some $\theta\in (0, 1]$.
The result is regarded as a nontrivial generalization of \cite{T24} to the non-autonomous case although the weight function is restricted to the isotropic one.
Indeed, for the case $\alpha\in [0, 1)$, the rate of decay/growth is the same as the one for the Oseen semigroup obtained in \cite[Theorem 2.5]{T24}.
\par
The non-autonomous character makes the problem quite involved, although the strategy of the proof is similar to that of \cite{Hi20, HiShi, I, KK15, KShi, T24}, in which the procedure consists of deduction of two estimates: One is the estimate near the boundary $\partial D$ (so-called the local energy decay estimate), the other is the estimate near spatial infinity.
We observe that the data taken from the weighted spaces possess better summability at infinity, which yields a local energy decay estimate with faster rate than the known one (\cite[Proposition 6.2]{Hi20}), see Proposition \ref{wei-led} below.
This enables us to show the gradient estimate of the desired form
\begin{align}\label{intro-grad}
\|\rho^{\alpha}\nabla T(t, s)f\|_q\le C(t-s)^{-\frac{1}{2}}\|\rho^{\alpha}f\|_q+C(t-s)^{-\frac{1}{2}}(1+t-s)^{\alpha}\|f\|_q
\end{align}
for $\alpha\in [0, 1)$ and $3/(3-\alpha)<q\le 3/(1-\alpha)$, from which we see that the upper bound $3/(1-\alpha)$ is relaxed when $\alpha\in (0, 1)$, see the item 1 of Remark \ref{rmk-thm1-1}.
This finding is due to Tomoki Takahashi \cite{T24} for the Oseen semigroup with anisotropic weights.
The restriction $\alpha<1$ arises from the boundedness of the Fujita-Kato projection associated with the Helmholtz decomposition: The inequality
\begin{align}\label{intro-FK}
\|\rho^{\alpha}P_{G}f\|_{L^q(G)}\le C\|\rho^{\alpha}f\|_{L^q(G)}
\end{align}
holds true if and only if $3/(3-\alpha)<q<\infty$, where $G$ is either the whole space $\mathbb R^3$ or the exterior domain $D$ under consideration, see also \eqref{FK-wei}.
We use \eqref{intro-FK} with $q<3/2$ and $G=\mathbb R^3$ to deduce the desired estimates near spatial infinity, which is possible only if $\alpha<1$.
By utilizing the local energy decay estimate with faster rate, we also prove the weighted estimates for the case $\alpha\in [1, 3)$, which would be less useful, nevertheless, they are of independent interest.\par
Although we focus on the translational case, we also prove the weighted estimates of even more satisfactory form when the rigid motion is purely rotatory.
Indeed, the growth $(1+t-s)^{\alpha}$ in \eqref{intro-grad} which stems from the translation does not appear, see Remark \ref{rmk-th1-3}.
As for the Stokes semigroup, where the rotation is also absent, the decay estimates with the isotropic weight were found first by Kobayashi and Kubo \cite{KK15}, see also Farwig and Tsuda \cite{FT26}.
Our result mentioned in Remark \ref{rmk-th1-3} can be regarded as a generalization of their results to the non-autonomous situation.
In the nonlinear context, however, these fine estimates will not play any role in this paper because we consider the time-periodic flow $v$ satisfying \eqref{subcrit} with $\lambda>0$, where $\lambda$ is given by \eqref{def-lambda}.
In fact, if the translation is absent, the basic flow is of scale-critical class such as $v(t)\in L^{3, \infty}(D)$ (weak-$L^3$ space), which cannot be handled by the method developed in the final section; indeed, we should take a similar way to the one found by Yamazaki \cite{Y}.
\par
Our linear theory for the case $\alpha\in[0, 1)$ can be applied to showing the global well-posedness of the nonlinear problem \eqref{ns}.
Let $(v, p_v)$ be the time-periodic solution fulfilling \eqref{subcrit} with $\lambda>0$ obtained in \cite{Ga22, T26}, and set
\begin{align*}
w:=u-v, \qquad p:=p_u-p_v,
\end{align*}
which obey the following problem
\begin{align}\label{int-ns-dist}
\begin{split}
&\left.
\begin{array}{rl}
\partial_tw+w\cdot\nabla w&=\Delta w +(\eta+\omega\times x)\cdot\nabla w-\omega\times w
-v\cdot\nabla w-w\cdot\nabla v-\nabla p \\
\mbox{div $w$}&=0 \\
\end{array}
\right\} \;\;\mbox{in $D\times (0,\infty)$},  \\
&\quad w|_{\partial D}=0, \qquad
\lim_{|x|\to\infty}w=0,\qquad
w(\cdot, 0)=w_0:=u_0-v(\cdot, 0).
\end{split}
\end{align}The second aim of this paper is to construct a solution to \eqref{int-ns-dist} globally in time within the weighted framework. 
In Theorem \ref{thm2} below, we show that the problem \eqref{int-ns-dist} admits a unique solution $w$ which enjoys
\begin{align}
&\|\rho^{\alpha}w(t)\|_{q}=o\left(t^{-\frac{1}{2}+\frac{3}{2q}+\alpha}\right),\qquad q\in [3, \infty],\label{int-dist-beh}\\
&\|\rho^{\alpha}\nabla w(t)\|_{r}=o\left(t^{-\frac{1}{2}+\alpha}\right),\qquad\;\; r\in [3, \infty),\label{int-inf-beh}
\end{align}
as $t\to\infty$ under suitable smallness of the initial disturbance $w_0$ as well as the periodic flow $v$, where $\alpha\in[0, 1)$.
If, in particular, $\alpha\in [0, 1/2)$, then our result tells us that the periodic solution is asymptotically stable with respect to $\|\rho^{\alpha}(\cdot)\|_{\infty}$.
Moreover, as a by-product, the spatial-temporal behavior of the disturbance is also available.
The behaviors \eqref{int-dist-beh}--\eqref{int-inf-beh} with $\alpha=0$ recover the results in \cite{T26}.
\par
Unfortunately, it is still out of reach to obtain a global strong solution to \eqref{int-ns-dist}.
First of all, the case $\alpha=1$ is missing for the linear estimate of the desired form.
Moreover, there is another difficulty. 
Since the basic flow $v$ is of scale-subcritical class, we treat the linear terms $P_D(v\cdot\nabla w+w\cdot\nabla v)$ as perturbation from the evolution operator as in \cite{ES05, T24}. 
To show appropriate estimates of the corresponding Duhamel term, we have to estimate $\|\rho^{\alpha}P_D(v\cdot \nabla w+w\cdot\nabla v)\|_q$ with $q<3/2$, but the boundedness \eqref{intro-FK} of the projection in the weighted spaces needs the restriction $\alpha<1$, see \eqref{def-delta} and Lemma \ref{stab-lem-duha}.
Thus, it seems that somewhat different approach would be needed to analyze the full nonlinear problem even if one can show the $L^q$-$L^r$ estimates of the form \eqref{dec-0-alpha1}--\eqref{dec-1-alpha1} for $\alpha\ge 1$.
\par
This paper is organized as follows.
In section \ref{sect-result}, after introducing some basic notation, we give our main results.
Section \ref{sect-prelim} is devoted to collecting some auxiliary lemmas.
In section \ref{sect-led}, we prove the decay estimates near the obstacle (the local energy decay estimates) with weighted data.
The weighted $L^q$-$L^r$ estimates globally in time are proved in section \ref{sect-dec}.
In the final section, we prove the main theorem on the large time behavior of the Navier-Stokes flow past a rigid body.

\section{Main Results}\label{sect-result}
We first introduce basic notation.
For $q\in [1, \infty], k\in \mathbb Z_{\ge0}$ and a domain $G\subset \mathbb R^3$, we denote the standard Lebesgue and Sobolev spaces by $L^q(G)$ and by $W^{k, q}(G)$.
We abbreviate the $L^q$-norm $\|\cdot\|_{q, G}=\|\cdot\|_{L^q(G)}$ and even $\|\cdot\|_q=\|\cdot\|_{L^q(D)}$, where $D$ is the exterior domain under consideration with $C^{1, 1}$-boundary $\partial D$.
Without loss of generality, we always assume that
\begin{align*}
\mathbb R^3\setminus D\subset B_1,
\end{align*}
where $B_R$ denotes the open ball centered at the origin with radius $R>0$.
We set $D_R=D\cap B_R$ for $R\ge 1$. 
The class $C_0^\infty(G)$ comprises all $C^\infty$-functions with compact support in $G$ and the space $W^{k,q}_0(G)$
with integer $k\geq 0$ stands for the completion of $C_0^\infty(G)$ in $W^{k,q}(G)$. 
By $W^{-1,q}(G)$ we denote the dual space of $W^{1,q^\prime}_0(G)$, where $1/q^\prime+1/q=1$.
For $q\in (1, \infty)$ and $r\in [1, \infty]$, the Lorentz space $L^{q, r}(G)$ is the set of measurable functions $f$ satisfying $\|f\|_{q, r}<\infty,$ where
\begin{align*}
\|f\|_{q, r}:=
\begin{cases}
\displaystyle \left(\int_{0}^{\infty}\left(t^{\frac{1}{q}}f^{*}(t)\right)^{r}\,\frac{dt}{t}\right)^{\frac{1}{r}},\qquad &1\le r<\infty,\\[5mm]
\displaystyle \sup_{t>0}t^{\frac{1}{q}}f^{*}(t),&r=\infty.
\end{cases}
\end{align*}
Here, the decreasing rearrangement is given by
\begin{align*}
&f^{*}(t):=\inf\{\sigma>0;\;  m(\sigma, f)\le t\},\quad t>0,\\
&m(\sigma, f):=\mu\left(\{x\in G;\;  |f(x)|>\sigma\}\right),\quad \sigma>0,
\end{align*}
with $\mu$ being the Lebesgue measure on $\mathbb R^3$, see \cite{BL, G, L} for details.\par
Let $X$ be a reflexive Banach space and $\langle\cdot, \cdot\rangle$ a duality pairing between $X$ and its dual $X^*$.
For interval $I\subset \mathbb R$, we define $C_w(I; X):=\{u: I\to X;\; \langle u(\cdot), \psi\rangle\in C(I)\; \forall \psi\in X^{*}\}$, which is the set of $X$-valued weakly continuous functions on $I$. 
If $X$ is not reflexive, the space $C_w(I; X)$ should be replaced by $C_{w^*}(I; X)$, where $\langle\cdot, \cdot\rangle$ stands for the pairing between $X$ and its pre-dual.
By $BC(I; X)$, we denote the space of bounded continuous functions on $I$ with values in $X$.\par
We adopt the same symbols for denoting scalar and vector function spaces as long as there is no confusion.
Let $X_1$ and $X_2$ be two Banach spaces.
Then ${\mathcal L}(X_1,X_2)$ stands for the Banach space consisting of all bounded linear operators from $X_1$ into $X_2$.
We simply write ${\mathcal L}(X_1)={\mathcal L}(X_1,X_1)$.
Several positive constants are denoted by $C$, which may change from line to line.

Let us introduce the solenoidal function spaces.
Let $G\subset \mathbb R^3$ be one of the following domains;
the exterior $C^{1,1}$-domain $D$ under consideration, a bounded domain $D_R$ for $R\geq 1$, and the whole space $\mathbb R^3$.
The class $C_{0,\sigma}^\infty(G)$ consists of all divergence-free vector fields being in $C_0^\infty(G)$.
Let $1<q<\infty$, and denote by $L^q_\sigma(G)$ the completion of $C_{0,\sigma}^\infty(G)$, then it is characterized as
\begin{equation}
L^q_\sigma(G)=\{u\in L^q(G);\; \mbox{div $u$}=0,\, \nu\cdot u|_{\partial G}=0\},
\label{sole-chara}
\end{equation}
where $\nu$ stands for the outer unit normal to $\partial G$ and $\nu\cdot u$ is understood in the sense of normal trace
on $\partial G$ (this boundary condition is absent when $G=\mathbb R^3$).
We also define the solenoidal Lorentz spaces by
\begin{align}\label{sole-lor}
L^{q, r}_{\sigma}(G):=\left(L^{q_0}_{\sigma}(G), L^{q_1}_{\sigma}(G)\right)_{\theta, r}
\end{align}
with
\begin{align*}
1<q_0<q<q_1<\infty,\quad \frac{1}{q}=\frac{1-\theta}{q_0}+\frac{\theta}{q_1},\quad 1\le r\le \infty,
\end{align*}
where $(\cdot, \cdot)_{\theta, r}$ stands for the real interpolation functor.
The right-hand side of \eqref{sole-lor} is independent of choice of $\{q_0, q_1\}$, so that the space $L^{q, r}_{\sigma}(G)$ is well-defined.
The space of $L^q$-vector fields admits the Helmholtz decomposition
\begin{equation}
L^q(G)=L^q_\sigma(G)\oplus \{\nabla p\in L^q(G);\; p\in L^q_{\rm loc}(\overline{G})\}
\label{helm}
\end{equation}
which was proved by Fujiwara and Morimoto \cite{FM}, Miyakawa \cite{Mi}, and Simader and Sohr \cite{SiSo}.
By $P_G=P_{G,q}: L^q(G)\to L^q_\sigma(G)$, we denote the Fujita-Kato bounded projection associated with the decomposition \eqref{helm}.
Note the duality relation $(P_{G,q})^*=P_{G,q^\prime}$ as well as $L^q_\sigma(G)^*=L^{q^\prime}_\sigma(G)$, where
$1/q^\prime+1/q=1$.
We simply write $P=P_D$ for the exterior domain $D$ under consideration.
When $G=\mathbb R^3$, we see that $P_{\mathbb R^3}=I+{\mathcal R}\otimes {\mathcal R}$, where
${\mathcal R}=\nabla (-\Delta)^{-1/2}$ denotes the Riesz transform and $I$ stands for the identity operator.
\par
We now introduce the weighted spaces.
Let $G=\mathbb R^3$ or the exterior domain $D$ under consideration.
For $q\in [1, \infty]$ and $\alpha\ge 0$, the weighted Lebesgue space $L^q_{\alpha}(G)$ consists of all measurable functions $f$ satisfying $\rho^{\alpha}f\in L^q(G)$, where $\rho$ stands for the weight function given by \eqref{weight}.
We also define the weighted solenoidal vector fields $L^{q}_{\alpha, \sigma}(G)$ by
\begin{align*}
L^{q}_{\alpha, \sigma}(G):=L^{q}_{\alpha}(G)\cap L^{q}_{\sigma}(G)=
\begin{cases}
\{u\in L^q_{\alpha}(\mathbb R^3);\; \mathrm{div}\;u=0\},\qquad &G=\mathbb R^3,\\
\{u\in L^q_{\alpha}(D);\; \mathrm{div}\;u=0, \nu\cdot u|_{\partial D}=0\},\qquad &G=D,\\
\end{cases}
\end{align*}
Both $L^{q}_{\alpha}(G)$ and $L^q_{\alpha, \sigma}(G)$ are Banach space with norm $\|\rho^{\alpha}(\cdot)\|_{q, G}$.
Assuming, in particular, that $\alpha\in (0, 3)$ and $q\in (\frac{3}{3-\alpha}, \infty)$, then the weight function $\rho^{\alpha q}$ belongs to the Muckenhoupt class $\mathcal A_q(G)$, see Farwig and Sohr \cite[Section 2]{FS97} for the details.
It should be noted that the boundedness of the Fujita-Kato projection $P_G$
\begin{align}\label{FK-wei}
\|\rho^{\alpha}P_{G}f\|_{q, G}\le C\|\rho^{\alpha}f\|_{q, G},\qquad G=\mathbb R^3, D
\end{align}
holds true if and only if $\alpha\in (0, 3)$ and $q\in (\frac{3}{3-\alpha}, \infty)$.
Specifically, for the case $G=\mathbb R^3$, the boundedness \eqref{FK-wei} can also be understood in terms of that of the Riesz transform $\mathcal R=\nabla(-\Delta)^{-\frac{1}{2}}$; the validity of the inequality $\|\rho^{\alpha}\mathcal Rf\|_{q, \mathbb R^3}\le C\|\rho^{\alpha}f\|_q$ is equivalent to the aforementioned condition on $\alpha$ and $q$, see Stein  \cite[Chapter V]{St} and Torchinsky \cite[Chapter IX]{Tor}.
For such $\alpha$ and $q$, the authors of \cite{FS97} showed the Helmholtz decomposition of the weighted space
\begin{align*}
L^{q}_{\alpha}(G)=L^{q}_{\alpha, \sigma}(G)\oplus \{\nabla p\in L^q_{\alpha}(G);\; p\in L^q_{\mathrm{loc}}(\overline{G})\}
\end{align*}
and denseness of $C^{\infty}_{0, \sigma}(G)$ in $L^{q}_{\alpha, \sigma}(G)$.
When $G=D$, the decomposition and denseness for the specific weight $\rho^{\alpha}$ were already shown by Specovius-Neugebauer \cite[Corollary 3, Lemma 9]{SN90}, even before \cite{FS97}.
For later use, we define the set of exponent as
\begin{align}
\mathcal{E}(\alpha):=&\left\{(q, r);\;\frac{3}{3-\alpha}<q\le r\le \infty, q\neq \infty \right\},\label{exp-qr}
\end{align}
where $\alpha\in [0, 3)$.

We introduce the generator associated with the non-autonomous linearized system \eqref{linear}.
Suppose that the rigid motion $\eta(t)+\omega(t)\times x$ fulfills \eqref{ass1-rigid} for some $\theta\in (0,1]$ as in \cite{HR14}. 
Let $1<q<\infty$. We define the family of linear operators $\{L(t)\}_{t\geq 0}$ on $L^q_\sigma(D)$ by
\begin{equation}
\begin{split}
D_q(L(t))&=\{u\in D_q(A);\; (\omega(t)\times x)\cdot\nabla u\in L^q(D)\}, \\
L(t)u&=P{\mathcal L}(t)u 
=Au-(\eta(t)+\omega(t)\times x)\cdot\nabla u+\omega(t)\times u,
\end{split}
\label{ig}
\end{equation}
where ${\mathcal L}(t)$ is the differential operator given by \eqref{pre-proj-0} and
$A$ denotes the Stokes operator on $L^q_\sigma(D)$, that is,
\begin{equation*}
D_q(A)=L^q_\sigma(D)\cap W^{1,q}_0(D)\cap W^{2,q}(D), \qquad Au=-P\Delta u.
\end{equation*}
In \eqref{ig} we have taken into account
\begin{equation*}
P\big[(\eta+\omega\times x)\cdot\nabla u-\omega\times u\big]=(\eta+\omega\times x)\cdot\nabla u-\omega\times u
\end{equation*}
and this is because 
\begin{equation}
\begin{split}
&\mbox{div $\big[(\eta+\omega\times x)\cdot\nabla u-\omega\times u\big]$}=(\eta+\omega\times x)\cdot\nabla\mbox{div $u$}=0, \\
&\nu\cdot\big[(\eta+\omega\times x)\cdot\nabla u-\omega\times u\big]\big|_{\partial D}=0,
\end{split}
\label{grund}
\end{equation}
see \eqref{sole-chara} and \eqref{helm}.
The latter of \eqref{grund} in our context was observed by Galdi and Silvestre \cite{GS05} when $q=2$.
Their finding is still valid even for $1<q<\infty$ since
the normal trace of the drift term vanishes, namely,
$\nu\cdot (\partial_ju)|_{\partial D}=0$ for every $j\in \{1,2,3\}$ and $u$ being in $L^q_\sigma(D)\cap W^{1,q}_0(D)$
that coincides with 
the completion of $C_{0,\sigma}^\infty(D)$ in $W^{1,q}(D)$.
As for the coincidence of those spaces, which the point here, see, for instance,
Galdi \cite[Theorem III.4.2]{Ga-b}; further references are found in \cite[Section III.7]{Ga-b}.
The initial value problem \eqref{linear} is now formulated as
\begin{equation*}
\frac{du}{dt}+L(t)u=0, \quad t\in (s,\infty); \qquad u(s)=f
\end{equation*}
in $L^q_\sigma(D)$.
\par
Hansel and Rhandi \cite{HR14} proved the generation of the evolution operator $T(t, s)$ on $L^q_{\sigma}(D)$ and the $L^q$-$L^r$ smoothing estimate near the initial time $s$.
Later on, further regularity properties with the aid of the weighted estimates were established by Hishida and the present author \cite{AH25}. 
\begin{proposition}[\cite{AH25, HR14}]\label{prop1-gen-sm}
Suppose that $\eta$ and $\omega$ fulfill \eqref{ass1-rigid} for some $\theta\in (0, 1]$. Let $1<q<\infty$. The operator family $\{L(t)\}_{t\ge0}$ generates an evolution operator $\{T(t, s)\}_{t\ge s\ge0}$ on $L^q_{\sigma}(D)$ such that $T(t, s)$ is a bounded operator from $L^{q}_{\sigma}(D)$ into itself with the semigroup property
\begin{align*}
T(t, \tau)T(\tau, s)=T(t, s)\quad (t\ge\tau\ge s\ge0) ; \qquad T(s, s)=I,
\end{align*}
in $\mathcal L(L^{q}_{\sigma}(D))$ with $I$ being the identity operator and that the map
\begin{align*}
\{t\ge s\ge 0\}\ni (t, s)\mapsto T(t, s)f\in L^q_{\sigma}(D)
\end{align*}
is continuous for every $f\in L^q_{\sigma}(D)$. Furthermore, we have the following properties.
\begin{enumerate}
\item Let $\alpha\in [0, 3), (q, r)\in \mathcal{E}(\alpha)$ and $j=0, 1$. For each $\mathcal T\in (0, \infty)$ and $m\in (0, \infty)$, there is a constant $C=C(\mathcal T, m, \alpha, q, r, \theta, D)>0$ such that
\begin{align}\label{wei-est1}
\|\rho^{\alpha}\nabla^{j}T(t, s)f\|_{r}\le C(t-s)^{-\frac{j}{2}-\frac{3}{2}\left(\frac{1}{q}-\frac{1}{r}\right)}\|\rho^{\alpha}f\|_{q}
\end{align}
for all $(t, s)$ with $0\le s<t\le \mathcal T$ and $f\in L^q_{\alpha, \sigma}(D)$ whenever 
\begin{align}\label{ever-loc}
\sup_{0\le t\le\mathcal T}(|\eta(t)|+|\omega(t)|)+\sup_{0\le s<t\le \mathcal T}\frac{|\eta(t)-\eta(s)|+|\omega(t)-\omega(s)|}{(t-s)^{\theta}}\le m
\end{align}
is satisfied, where the weight function $\rho$ is given by \eqref{weight}. Moreover, we have
\begin{align}\label{wei-small1}
\lim_{t\to s}\, (t-s)^{\frac{j}{2}+\frac{3}{2}\left(\frac{1}{q}-\frac{1}{r}\right)}\|\rho^\alpha\nabla^jT(t,s)f\|_r=0
\end{align}
for every $f\in L^q_{\alpha,\sigma}(D)$ except the case $(j,r)=(0,q)$;
in that case, it holds that
\begin{align}
\lim_{t\to s}\|\rho^\alpha\bigl(T(t,s)f-f\bigr)\|_q=0
\label{IC-1}
\end{align}
for every $f\in L^q_{\alpha,\sigma}(D)$.
\item Let $3/2<q<\infty$. 
For every $f\in L^q_{1,\sigma}(D)$, the function 
$T(t,s)f$ 
is of class 
\begin{align}\label{st-cl}
T(\cdot, s)f\in C^1((s, \infty); L^q_{\sigma}(D))\cap C((s, \infty); W^{2, q}(D)),
\end{align}
with
\begin{align}\label{str}
\partial_tT(t, s)f+L(t)T(t, s)f=0,\qquad t\in (s, \infty)
\end{align}
in $L^q_{\sigma}(D)$.
For each ${\mathcal T}\in (0,\infty)$ and $m\in (0,\infty)$, there is a constant $C=C({\mathcal T},m,q,\theta,D)>0$ such that 
\begin{align*}
\|\partial_tT(t,s)f\|_q+\|L(t)T(t,s)f\|_{q}+\|\nabla p(t)\|_q  
\leq C(t-s)^{-1}\|f\|_q+C(t-s)^{-1/2}\|\rho f\|_q
\end{align*}
for all $(t,s)$ with $0\leq s<t\le {\mathcal T}$ and $f\in L^q_{1,\sigma}(D)$ whenever \eqref{ever-loc} is satisfied, where $p(t)$ is the pressure associated with $T(t,s)f$.
\end{enumerate}
\end{proposition}

\begin{remark}
We regard $T(t, s)$ as the evolution operator on $L^q_{\sigma}(D)$, not on the weighted space $L^{q}_{\alpha, \sigma}(D)$ at all in the present paper.
We do not know whether the continuity 
\begin{align*}
\{t\ge s\ge 0\}\ni (t, s)\mapsto T(t, s)f\in L^{q}_{\alpha, \sigma}(D)
\end{align*}
holds, where $\alpha\in(0, 3), q\in (\frac{3}{3-\alpha}, \infty)$ and $f\in L^q_{\alpha, \sigma}(D)$.
\end{remark}

We next consider the $L^q$-$L^r$ decay estimates of the evolution operator $T(t, s)$.
Under the condition \eqref{ass2-rigid} for some $\theta\in (0, 1]$, we set
\begin{align}\label{norm-rigid}
\begin{split}
&|(\eta, \omega)|_{\infty}:=\sup_{t\ge0}(|\eta(t)|+|\omega(t)|),\\
&|(\eta, \omega)|_{\theta}:=\sup_{t>s\ge 0}\frac{|\eta(t)-\eta(s)|+|\omega(t)-\omega(s)|}{(t-s)^{\theta}}.
\end{split}
\end{align}
The next proposition was proved by Hishida \cite{Hi18, Hi20}.
\begin{proposition}[\cite{Hi18, Hi20}]\label{prop2-dec}
Suppose that $\eta$ and $\omega$ fulfill \eqref{ass2-rigid} for some $\theta\in (0, 1]$.
Let $1<q<\infty$ and $r\in [q, \infty]$, except $r\neq \infty$ for \eqref{lqlr-hi-1} below. 
For each $m\in (0, \infty)$, there is a constant $C=C(m, q, r, \theta, D)>0$ such that
\begin{align}
\|T(t, s)f\|_r&\le C(t-s)^{-\frac{3}{2}\left(\frac{1}{q}-\frac{1}{r}\right)}\|f\|_q\label{lqlr-hi-0}\\
\|\nabla T(t, s)f\|_r&\le C(t-s)^{-\frac{1}{2}-\frac{3}{2}\left(\frac{1}{q}-\frac{1}{r}\right)}(1+t-s)^{\max\left\{\frac{1}{2}-\frac{3}{2r}, 0\right\}}\|f\|_q\label{lqlr-hi-1}
\end{align}
for all $t>s\ge0$ and $f\in L^q_{\sigma}(D)$, whenever
\begin{align}\label{ever}
|(\eta, \omega)|_{\infty}+|(\eta, \omega)|_{\theta}\le m
\end{align}
is satisfied.
\end{proposition}

We are in a position to give our first theorem on the weighted estimates of the evolution operator $T(t, s)$ globally in time.
The main assertions are \eqref{dec-0-alpha1}--\eqref{dec-1-alpha1} below for the case $\alpha\in [0, 1)$;
in fact, they are of the desired form in which the term $(1+t-s)^{\alpha}$ appears only with $\|f\|_q$.
As clarified by Tomoki Takahashi \cite{T24}, who discussed the large time behavior of the Oseen semigroup with anisotropic weights $(1+|x|)^{\alpha_1}(1+|x|-x_1)^{\alpha_2}$ reflecting the wake structure, estimates of the form \eqref{dec-0-alpha1}--\eqref{dec-1-alpha1} are useful in analyzing the nonlinear problem.
Indeed, our estimates \eqref{dec-0-alpha1}--\eqref{dec-1-alpha1} recover his corresponding results with $\alpha_2=0$ and are utilized in the proof of Theorem \ref{thm2}.
If in particular $\eta=0$, even better preferable estimates are available, see Remark \ref{rmk-th1-3} below, however, we need to study the translational case as mentioned in the introductory section.
When $\alpha\in [1, 3)$, unfortunately, it seems to be difficult to avoid even more growth than $(1+t-s)^{\alpha}$ in the second term of \eqref{dec-0-alpha1}--\eqref{dec-1-alpha1}, nevertheless, estimates of the form \eqref{dec-0-alpha2}--\eqref{dec-1-alpha2} below can be deduced.
Although those latter estimates for $\alpha\ge 1$ would not be very useful, they might be of independent interest and, thus, are included in the following theorem.

\begin{theorem}\label{thm1}
Suppose that $\eta$ and $\omega$ fulfill \eqref{ass2-rigid} for some $\theta\in (0, 1]$. Let $\alpha\in [0, 3)$ and $(q, r)\in\mathcal{E}(\alpha)$, see \eqref{exp-qr}.
\begin{enumerate}
\item  Let $\alpha\in [0, 1)$.  For each $m\in (0, \infty)$, there is a constant $C=C(m, \alpha, q, r, \theta, D)>0$ such that 
\begin{align}\label{dec-0-alpha1}
\|\rho^{\alpha}T(t, s)f\|_r\le C(t-s)^{-\frac{3}{2}\left(\frac{1}{q}-\frac{1}{r}\right)}\|\rho^{\alpha}f\|_q+(t-s)^{-\frac{3}{2}\left(\frac{1}{q}-\frac{1}{r}\right)}(1+t-s)^{\alpha}\|f\|_q
\end{align}
for all $t>s\ge 0$ and $f\in L^q_{\alpha, \sigma}(D)$ whenever \eqref{ever} is satisfied.
\item Let $\alpha\in [0, 1)$ and $r\neq \infty$. For each $m\in (0, \infty)$, there is a constant $C=C(m, \alpha, q, r, \theta, D)>0$ such that 
\begin{align}\label{dec-1-alpha1}
\|\rho^{\alpha}\nabla T(t, s)f\|_r
&\le C(t-s)^{-\frac{1}{2}-\frac{3}{2}\left(\frac{1}{q}-\frac{1}{r}\right)}(1+t-s)^{\max\left\{0, \frac{1}{2}-\frac{3}{2r}-\frac{\alpha}{2}\right\}}\Bigl(\|\rho^{\alpha}f\|_q+(1+t-s)^{\alpha}\|f\|_{q}\Bigr)
\end{align}
for all $t>s\ge 0$ and $f\in L^q_{\alpha, \sigma}(D)$ whenever \eqref{ever} is satisfied.
\item Let $\alpha\in [1, 3)$.
For each $m\in (0, \infty)$, there is a constant $C=C(m, \alpha, q, r, \theta, D)>0$ such that 
\begin{align}\label{dec-0-alpha2}
\|\rho^{\alpha}T(t, s)f\|_r\le C(t-s)^{-\frac{3}{2}\left(\frac{1}{q}-\frac{1}{r}\right)}(1+t-s)^{\alpha}\|\rho^{\alpha}f\|_q
\end{align}
for all $t>s\ge 0$ and $f\in L^q_{\alpha, \sigma}(D)$ whenever \eqref{ever} is satisfied.
\item Let $\alpha\in [1, 3)$ and $r\neq\infty$.
Assume, in addition, $1/q-1/r<1/3$ when $\alpha\in [1, 2)$.
For each $m\in (0, \infty)$, there is a constant $C=C(m, \alpha, q, r, \theta, D)>0$ such that 
\begin{align}\label{dec-1-alpha2}
\|\rho^{\alpha}\nabla T(t, s)f\|_r\le C(t-s)^{-\frac{1}{2}-\frac{3}{2}\left(\frac{1}{q}-\frac{1}{r}\right)}(1+t-s)^{\alpha}\|\rho^{\alpha}f\|_q
\end{align}
for all $t>s\ge 0$ and $f\in L^q_{\alpha, \sigma}(D)$ whenever \eqref{ever} is satisfied.
\end{enumerate}
\end{theorem}

\begin{remark}\label{rmk-thm1-1}
We give two comments on \eqref{dec-0-alpha1} with $r=q$ for the case $\alpha\in (0, 1)$.
\begin{enumerate}
\item As a particular case of \eqref{dec-0-alpha1}, we have
\begin{align}\label{rmk-thm1-lqlq}
\|\rho^{\alpha}\nabla T(t, s)f\|_q
\le
(t-s)^{-\frac{1}{2}}\|\rho^{\alpha}f\|_q+(t-s)^{-\frac{1}{2}+\alpha}\|f\|_q,\quad \; \frac{3}{3-\alpha}<q\le \frac{3}{1-\alpha},
\end{align}
for all $t-s>2$ as well as $0\le s<t$ and $f\in L^q_{\alpha, \sigma}(D)$.
The corresponding estimate for the case $\alpha=0$ was proved by Hishida \cite{Hi20}, see Proposition \ref{prop2-dec}.
His result covers the specific cases studied by the following literature; \cite{I, MaSol} (the Stokes semigroup, both $\eta$ and $\omega$ are absent),  \cite{KShi} (the Oseen semigroup, $\eta$ is a nonzero constant and $\omega$ is absent), \cite{HiShi}  (the Stokes semigroup with rotating effect, $\eta$ is absent and  $\omega$ is a nonzero constant), \cite{Shi08} (the Oseen semigroup with rotating effect, $\eta$ and $\omega$ are nonzero constants which are parallel to each other).
For the Stokes semigroup, Maremonti and Solonnikov \cite{MaSol} proved that one can not avoid the restriction $q\le 3=n$ (space dimension) for the desired decay rate $(t-s)^{-\frac{1}{2}}$ in \eqref{rmk-thm1-lqlq} with $\alpha=0$. 
See also \cite{Hi11} for an alternative proof. 
With the aid of the weighted analysis, we can extend the upper bound for $q$ to obtain the desired rate when $0\le\alpha<1$.
This finding is due to Tomoki Takahashi \cite{T24} for the Oseen semigroup.
\item For the case $\alpha\in (0, 1)$, we can actually obtain even better rate of decay. 
We conclude \eqref{dec-1-alpha1} by combining the local energy decay estimate \eqref{led-wei} with the estimate \eqref{dec-1-infty-01} near spatial infinity.
If we look carefully into details of the deduction of \eqref{dec-1-infty-01}, we find
\begin{align*}
\|\rho^{\alpha}\nabla T(t, s)f\|_{q}\le C(t-s)^{\max\left\{-\frac{1}{2}, -\frac{3}{2q}-\frac{\alpha}{2}\right\}}\|\rho^{\alpha}f\|_q+C(t-s)^{\max\left\{-\frac{1}{2}+\alpha, -\frac{3}{2q}\right\}}\|f\|_q
\end{align*}
for $t-s>2$, which improves \eqref{dec-1-alpha1} with $r=q>3/(1-\alpha)$.
However, this estimate would be less useful because the rate of decay cannot be kept when using the semigroup property.
See also Remark \ref{lqlr-rate-imp} regarding the same observation for $L^q$-$L^r$ estimates with $r>q$ as above.
\end{enumerate}
\end{remark}

\begin{remark}\label{rmk-th1-2}
For the gradient estimate \eqref{dec-1-alpha2}, the case $1/3\le 1/q-1/r<1-\alpha/3$ is missing when $\alpha\in [1, 2)$.
For this case, indeed, what one can obtain by using the semigroup property is the following with worse rate:
\begin{align*}
\|\rho^{\alpha}\nabla T(t, s)f\|_{r}\le C(t-s)^{-\frac{1}{2}-\frac{3}{2}\left(\frac{1}{q}-\frac{1}{r}\right)}(1+t-s)^{2\alpha}\|\rho^{\alpha}f\|_q.
\end{align*}
On the contrary, the semigroup property works well for \eqref{dec-0-alpha1}--\eqref{dec-1-alpha1} in the sense that we do not lose any rate when using this property.
\end{remark}

\begin{remark}\label{rmk-th1-3}
Consider the specific case $\eta\equiv0$, in which the rigid motion of the body is purely rotatory. 
Then we have the following improvement under the same assumptions as in Theorem \ref{thm1}:
\begin{align}
&\|\rho^{\alpha}T(t, s)f\|_r\le C(t-s)^{-\frac{3}{2}\left(\frac{1}{q}-\frac{1}{r}\right)}\|\rho^{\alpha}f\|_{q},\label{rot-0-imp}\\
&\|\rho^{\alpha}\nabla T(t, s)f\|_r\le C(t-s)^{-\frac{1}{2}-\frac{3}{2}\left(\frac{1}{q}-\frac{1}{r}\right)}(1+t-s)^{\max\left\{0, \frac{1}{2}-\frac{3}{2r}-\frac{\alpha}{2}\right\}}\|\rho^{\alpha}f\|_{q}.\label{rot-1-imp}
\end{align}
This is because the weighted estimates of the evolution operator $U(t, s)$ in the whole space $\mathbb R^3$ are better when $\eta\equiv0$, see Remark \ref{rem-wh-st-rot}.
For the Stokes semigroup, this was already discovered by Kobayashi and Kubo \cite[Theorem 2, (i)]{KK15}, Tomoki Takahashi \cite[(44)]{T24}, Farwig and Tsuda \cite[Theorem 3.1]{FT26}; thus, \eqref{rot-0-imp}--\eqref{rot-1-imp} may be regarded as a generalization of their results to the non-autonomous case.
We will give a brief comment on the proof of \eqref{rot-0-imp}--\eqref{rot-1-imp} in Remark \ref{rem-ext-st-rot}. 
Nevertheless, this improvement will not be used in the rest of the paper.
In fact, as we know from \eqref{subcrit} with \eqref{def-lambda}, the spatial pointwise behavior of time-periodic flows for the case $\eta=0$ is the scale-critical rate $O(|x|^{-1})$, which cannot be covered by Theorem \ref{thm2} below.
\end{remark}

As an application of Theorem \ref{thm1}, we consider the large time behavior of solutions to the Navier-Stokes system around a basic flow in weighted spaces.
Specifically, as the basic flow, we take a time-periodic solution $(v, p_v)$ to the problem \eqref{ns-per-int} constructed by Galdi \cite{Ga22} and Tomoki Takahashi \cite{T26}.
They considered the problem \eqref{ns-per-int} with a given external force, but we assume in the following proposition that it is absent for simplicity. 
\begin{proposition}[\cite{Ga22, T26}]\label{prop-galtak}
Assume that $\eta$ and $\omega$ are $\mathcal T$-periodic with a period $\mathcal T>0$.
Let $\eta, \omega\in W^{2, 2}(0, \mathcal T; \mathbb R^3)$ satisfy either \eqref{ass1-per} or \eqref{ass2-per}, and suppose \eqref{def-lambda} with $\lambda>0$. 
There exists a constant $\widetilde{\varepsilon}=\widetilde{\varepsilon}(D, \mathcal T)>0$ such that if
\begin{align}\label{per-as-small}
\|(\eta, \omega)\|_{W^{2, 2}(0, \mathcal T)}\le \widetilde{\varepsilon},
\end{align}
then problem \eqref{ns-per-int} has one and only one $\mathcal T$-periodic solution $(v, p_v)$, which enjoys the pointwise decay estimates 
\begin{align*}
\sup_{t\in\mathbb R}|\nabla^jv(x, t)|\le \frac{C}{(1+|x|)^{1+\frac{j}{2}}(1+\lambda|x|+\lambda x_1)^{1+\frac{j}{2}}}\|(\eta, \omega)\|_{W^{2, 2}(0, \mathcal T)},
\end{align*}
for $j=0, 1$ and $x\in D$ with some constant $C=C(D, \mathcal T)>0$.
\end{proposition}
Let $(v, p_v)$ be the periodic solution obtained in Proposition \ref{prop-galtak}.
We find a solution to \eqref{ns} of the form
\begin{align*}
u=v+w,\qquad p_u=p_v+p,
\end{align*}
then the disturbance $(w, p)$ should obey the problem \eqref{int-ns-dist}.
Our objective is to show the existence of a unique solution to \eqref{int-ns-dist} globally in time and deduce the decay/growth rate of the solution $w$ in weighted Lebesgue spaces.
Applying the Fujita-Kato projection to \eqref{int-ns-dist}, the problem is formulated as
\begin{align*}
\begin{split}
&\frac{dw}{dt}+L(t)w+P(w\cdot\nabla w+v\cdot\nabla w+w\cdot\nabla v)=0,\qquad t\in (0, \infty),\\
&w(0)=w_0,
\end{split}
\end{align*}
which is converted into the integral equation
\begin{align}\label{ns-inteq}
w(t)=T(t, 0)w_0-\int_{0}^{t}T(t, \tau)P(w\cdot\nabla w+v\cdot\nabla w+w\cdot\nabla v)(\tau)\,d\tau.
\end{align}
Since the rigid motion is $\mathcal T$-periodic and since $\eta, \omega\in W^{2, 2}(0, \mathcal T; \mathbb R^3)\hookrightarrow C^1([0, \mathcal T]; \mathbb R^3)$, we have \eqref{ass2-rigid} with $\theta=1$ subject to
\begin{align*}
|(\eta, \omega)|_{\infty}+|(\eta, \omega)|_{1}\le C\|(\eta, \omega)\|_{W^{2, 2}(0, \mathcal T)},
\end{align*}
see \eqref{norm-rigid}, that allows us to use Theorem \ref{thm1}.
Now, our main result reads as follows.
For the case $\alpha=0$, the corresponding result was already obtained in \cite[Theorem 2.9]{T26}.
For the meaning of the latter condition in \eqref{dist-as-small}, see \cite[Remark 2.10]{T26}. 
\begin{theorem}\label{thm2}
In addition to the assumptions of Proposition \ref{prop-galtak}, let $\alpha\in [0, 1)$.
There exist constants $\varepsilon_0=\varepsilon_0(\alpha, D, \mathcal T)>0$ and $\varepsilon_1=\varepsilon_1(\alpha, D, \mathcal T)>0$ such that if $w_0\in L^3_{\alpha, \sigma}(D), \eta, \omega\in W^{2, 2}(0, \mathcal T; \mathbb R^3)$ and $\lambda\in (0, 1]$ satisfy
\begin{align}\label{dist-as-small}
\|\rho^{\alpha}w_0\|_3\le \varepsilon_0,\qquad \lambda^{-\frac{3}{4}}\|(\eta, \omega)\|_{W^{2, 2}(0, \mathcal T)}\le\varepsilon_1
\end{align}
as well as \eqref{per-as-small}, then \eqref{ns-inteq} admits a solution
\begin{align}\label{sol-conti-l3}
w\in BC([0, \infty);L^3_{\sigma}(D))
\end{align}
which enjoys
\begin{align}
&t^{\frac{1}{2}-\frac{3}{2q}}w\in BC((0, \infty); L^q_{\sigma}(D)),\qquad q\in (3, \infty),\label{sol-conti-0th}\\
&t^{\frac{1}{2}}w\in BC_{w^*}((0, \infty); L^{\infty}(D)),\label{sol-conti-0-infty}\\
&t^{\frac{1}{2}}\nabla w\in BC_{w}((0, \infty); L^3(D)),\label{sol-conti-1st}\\
&t^{\frac{1}{2}-\frac{3}{2q}}(1+t)^{-\alpha}\rho^{\alpha}w\in L^{\infty}(0, \infty; L^q(D)),\;\; q\in (3, \infty];\quad t^{\frac{1}{2}}(1+t)^{-\alpha}\rho^{\alpha}\nabla w\in L^{\infty}(0, \infty; L^3(D)).\label{sol-wei}
\end{align}
The solution is unique within the class
\begin{align}\label{uni-tak}
\{w\in BC([0, \infty); L^3_{\sigma}(D));\; t^{\frac{1}{2}}w\in BC_{w^*}((0, \infty); L^{\infty}(D)), t^{\frac{1}{2}}\nabla w\in BC_w((0, \infty); L^3(D))\}.
\end{align}
\par
Moreover, there exist constants $\varepsilon_2=\varepsilon_2(\alpha, D, \mathcal T)\in (0, \varepsilon_0]$ and $\varepsilon_3=\varepsilon_3(\alpha, D, \mathcal T)\in (0, \varepsilon_{1}]$ such that if 
\begin{align*}
\|\rho^{\alpha}w_0\|_3\le \varepsilon_2,\qquad \lambda^{-\frac{3}{4}}\|(\eta, \omega)\|_{W^{2, 2}(0, \mathcal T)}\le\varepsilon_3,
\end{align*}
then the solution $w(t)$ obtained above satisfies
\begin{align}
&\|w(t)\|_{q}=o\left(t^{-\frac{1}{2}+\frac{3}{2q}}\right), &&\|\nabla w(t)\|_{r}=o\left(t^{-\frac{1}{2}}\right),\label{sol-0-sob}\\
&\|\rho^{\alpha}w(t)\|_q=o\left(t^{-\frac{1}{2}+\frac{3}{2q}+\alpha}\right), &&\|\rho^{\alpha}\nabla w(t)\|_r=o\left(t^{-\frac{1}{2}+\alpha}\right)\quad \label{sol-wei-sob}
\end{align}
as $t\to\infty$ for every $q\in [3, \infty]$ and $r\in [3, \infty)$.
\end{theorem}

\section{Preliminaries}\label{sect-prelim}

In this section, we collect some known results  and auxiliary lemmas for later use. 

\subsection{Bogovski\u\i\,  operator}
Let us consider the boundary value problem for the divergence equation
\begin{align}
\mbox{div $w$}=f\quad\mbox{in $G$}, \qquad w|_{\partial G}=0,
\label{bvp-div}
\end{align} 
where $G$ is a bounded domain in $\mathbb R^n$.
There are many solutions to \eqref{bvp-div} under the compatibility condition $\int_Gf\,dx=0$ and appropriate regularity of $f$. 
In \cite{Bog}, Bogovski\u\i\, found a particular solution to \eqref{bvp-div}, which is useful to recover the solenoidal condition in cut-off arguments due to some fine properties of his solution.
The result needed in this paper is summarized in the following lemma, see \cite[Theorem 1]{Bog}, \cite[Theorem 2.4 (a)--(c)]{BorSo90}, 
\cite[Theorem III.3.3]{Ga-b}, \cite[Theorem 2.5]{GHH06-2} and the references therein. 
\begin{lemma}\label{lem-bog}
Let $G\subset \mathbb R^n$, $n\geq 2$, be a bounded domain with Lipschitz boundary $\partial G$.
There exists a linear operator $\mathbb B_G: C_0^\infty(G)\to C_0^\infty(G)^n$ with the following properties:
For every $q\in (1,\infty)$ and integer $k\geq 0$, there is a constant $C=C(q,k,G)>0$ being invariant under dilation
of the domain $G$ such that
\begin{equation*}
\|\nabla^{k+1}\mathbb B_Gf\|_{q,G}\leq C\|\nabla^kf\|_{q,G}
\end{equation*}
and that $w=\mathbb B_Gf$ is a solution to \eqref{bvp-div} if $\int_Gf\,dx=0$.
The operator $\mathbb B_G$, that we call the Bogovski\u\i\, operator, 
extends uniquely to a bounded operator from $W_0^{k,q}(G)$ to $W_0^{k+1,q}(G)^n$.
Furthermore, it also extends uniquely to a bounded operator from $W^{1,q^\prime}(G)^*$ to $L^q(G)^n$, namely,
\begin{equation}
\|\mathbb B_Gf\|_{q,G}\leq C\|f\|_{W^{1,q^\prime}(G)^*}
\label{bog-est2}
\end{equation}
with some constant $C=C(q,G)>0$, where $1/q^\prime+1/q=1$.  
\end{lemma}

\subsection{Whole space flow, $\{U(t, s)\}_{t\ge s\ge 0}$}\label{Upre}
We consider the initial value propblem for the non-autonomous system
\begin{align}
\begin{split}
&\partial_tu+{\mathcal L}(t)u+\nabla p=0,\qquad\mbox{div $u$}=0 \quad\mbox{in $\mathbb R^3\times (s,\infty)$}, \\
&\lim_{|x|\to\infty}u=0, \qquad u(\cdot,s)=f,
\end{split}
\label{linear-wh}
\end{align}
where ${\mathcal L}(t)$ is the differential operator given by \eqref{pre-proj-0}.
Because of $\mbox{\eqref{grund}}_1$ we see that $\nabla p=0$ within the class $\nabla p\in L^q(\mathbb R^3)$.
Hence, under the compatibility condition $\mbox{div $f$}=0$, the solution formula to \eqref{linear-wh} is obtained from the heat semigroup
\[
\big(e^{t\Delta}f\big)(x)=(4\pi t)^{-3/2}\left(e^{-|\cdot|^2/4t}*f\right)(x)
\]
simply by transformation of variables as follows, where $*$ stands for convolution in spatial variable.
It is expicitly described as 
\begin{equation*}
u(x,t)=
\big(U(t,s)f\big)(x)=
\Phi(t,s)\left(e^{(t-s)\Delta}f\right)\left(\Phi(t,s)^\top\Big(x+\int_s^t\Phi(t,\tau)\eta(\tau)\,d\tau\Big)\right),
\end{equation*}
where $\Phi(t,s)\in \mathbb R^{3\times 3}$ stands for the evolution matrix, being an orthogonal matrix,
to the differential equation
$\frac{d\varphi}{dt}=-\omega(t)\times\varphi$.
\par
The estimates of $U(t, s)$ we need in this paper are summarized as follows. See also \cite[Proposition 4.1]{AH25}, in which they considered the local estimate with $0\le s< t\le \mathcal T$ for $\mathcal T\in (0, \infty)$, nevertheless, we know \eqref{wei-wh-all} below from \cite[(4.6)]{AH25}.
It is obvious that \eqref{wei-wh-all} implies \eqref{wei-wh-unif}.
The second assertion of the following lemma was already shown in \cite[Lemma 3.2]{Hi20}, where \eqref{ass2-rigid} was imposed, however, the H\"older continuity played no role there as seen from the proof.
\begin{lemma}\label{so-far-whole}
Let $\eta, \omega\in C([0, \infty); \mathbb R^3)$ and $\alpha\ge0$.
\begin{enumerate}
\item Assume that $q\in [1, \infty]$ and $r\in[q, \infty]$. Let $j\ge0$ be an integer. For each $m_0\in (0, \infty)$, there is a constant $C_0=C_0(q, r, \alpha, j)>0$ such that
\begin{align}\label{wei-wh-all}
\|\rho^{\alpha}\nabla^j U(t, s)f\|_{r, \mathbb R^3}\le C_0(t-s)^{-\frac{j}{2}-\frac{3}{2}\left(\frac{1}{q}-\frac{1}{r}\right)}\left[\|\rho^{\alpha}f\|_{q, \mathbb R^3}+(t-s)^{\frac{\alpha}{2}}\left\{1+m_0^{\alpha}(t-s)^{\frac{\alpha}{2}}\right\}\|f\|_{q, \mathbb R^3}\right]
\end{align}
for all $t>s\ge 0$ and $f\in L^q_{\alpha}(\mathbb R^3)$, whenever
\begin{align}\label{ever2}
\sup_{t\ge0}|\eta(t)|\le m_0,
\end{align}
is satisfied, where $\rho$ is the weight function given by \eqref{weight}. Moreover, for each $\mathcal T\in (0, \infty)$, there is a constant $C_1=C_1(q, r, \alpha, j, m_0, \mathcal T)>0$ such that
\begin{align}\label{wei-wh-unif}
\|\rho^{\alpha} \nabla^j U(t, s)f\|_{r, \mathbb R^3}\le C_1(t-s)^{-\frac{3}{2}\left(\frac{1}{q}-\frac{1}{r}\right)-\frac{j}{2}}\|\rho^{\alpha}f\|_{q, \mathbb R^3}
\end{align}
for all $(t, s)$ with $t-s\le \mathcal T$ as well as $0\le s<t$ and $f\in L^q_{\alpha}(\mathbb R^3)$, whenever \eqref{ever2} is satisfied.
\item Let $\eta, \omega\in C([0, \infty); \mathbb R^3)\cap L^{\infty}(0, \infty; \mathbb R^3)$ and $R>0$. For every $f\in L^{q}_{\sigma}(\mathbb R^3)$, we have $U(\cdot, s)f\in C^1((s, \infty); W^{-1, q}(B_R))$. For each $\mathcal T\in (0, \infty)$ and $m\in (0, \infty)$, there is a constant $C=C(\mathcal T, m,q, R)>0$ such that
\begin{align}
\|\partial_tU(t, s)f\|_{W^{-1, q}(B_R)}\le C(t-s)^{-\frac{1}{2}}\|f\|_{q, \mathbb R^3}\label{weak-wh-unif}
\end{align}
for all $(t, s)$ with $t-s\le \mathcal T$ as well as $0\le s<t$, whenever 
\begin{align*}
|(\eta, \omega)|_{\infty}\le m
\end{align*}
is satisfied, see \eqref{norm-rigid}.
\end{enumerate}
\end{lemma}

\begin{remark}\label{rem-wh-st-rot}
Let us consider the case $m_0=0$, that is, the translational velocity $\eta$ is absent, see \eqref{ever2}. 
Then \eqref{wei-wh-all} implies that 
\begin{align*}
\|\rho^{\alpha}\nabla^jU(t, s)f\|_{r, \mathbb R^3}\le C(t-s)^{-\frac{j}{2}-\frac{3}{2}\left(\frac{1}{q}-\frac{1}{r}\right)}\left[\|\rho^{\alpha}f\|_{q, \mathbb R^3}+(t-s)^{\frac{\alpha}{2}}\|f\|_{q, \mathbb R^3}\right].
\end{align*}
However, the term $(t-s)^{\frac{\alpha}{2}}$ can be removed with the aid of Lemma \ref{lem-low-int} below, when assuming $\alpha\in [0, 3)$ and $(q, r)\in\mathcal{E}(\alpha)$. 
In view of the proof of Lemma \ref{so-far-whole}, see \cite[Proposition 4.1]{AH25}, we have
\begin{align*}
\rho(x)^{\alpha}|\nabla^j \left(U(t, s)f\right)(x)|\le C(I+II)
\end{align*}
with
\begin{align*}
I:=(t-s)^{-\frac{j}{2}+\frac{\alpha}{2}}\left(\widetilde{G}_{t-s}*g\right)(x),\qquad II:=(t-s)^{-\frac{j}{2}}\left(\widetilde{H}_{t-s}*\rho^{\alpha}g\right)(x),
\end{align*}
where
\begin{gather*}
g(y):=f\left(\Phi(t, s)^{\top}y\right),\\
\widetilde{G}_{t-s}(x):=(t-s)^{-\frac{3}{2}}e^{-\frac{|x|^2}{16(t-s)}},\qquad \widetilde{H}_{t-s}(x):=(t-s)^{-\frac{3}{2}}e^{-\frac{|x|^2}{8(t-s)}},
\end{gather*}
and $*$ denotes the convolution on $\mathbb R^3$. 
It follows from Lemma \ref{lem-low-int} below that 
\begin{align*}
\|g\|_{\frac{3q}{3+\alpha q}, q, \mathbb R^3}=\|f\|_{\frac{3q}{3+\alpha q}, q, \mathbb R^3}\le C\|\rho^{\alpha} f\|_{q, \mathbb R^3}
\end{align*}
for $f\in L^{q}_{\alpha}(\mathbb R^3)$. 
Then, by the Young inequality on the Lorentz space (see \cite{ON62} for instance), we have
\begin{align*}
\|I\|_{r, \mathbb R^3}
=(t-s)^{-\frac{j}{2}+\frac{\alpha}{2}}\|\widetilde{G}_{t-s}*g\|_{r, \mathbb R^3}
&\le C(t-s)^{-\frac{j}{2}+\frac{\alpha}{2}}(t-s)^{-\frac{3}{2}\left(\frac{1}{q}-\frac{1}{r}\right)-\frac{\alpha}{2}}\|g\|_{\frac{3q}{3+\alpha q}, q, \mathbb R^3}\\
&\le C(t-s)^{-\frac{j}{2}-\frac{3}{2}\left(\frac{1}{q}-\frac{1}{r}\right)}\|\rho^{\alpha}f\|_{q, \mathbb R^3}.
\end{align*}
On the other hand, we see that $\|\rho^{\alpha}g\|_{q, \mathbb R^3}=\|\rho^{\alpha}f\|_{q, \mathbb R^3}$ due to $\eta=0$, yielding
\begin{align*}
\|II\|_{r, \mathbb R^3}=(t-s)^{-\frac{j}{2}}\|\widetilde{H}_{t-s}*\rho^{\alpha}g\|_{r, \mathbb R^3}
\le C(t-s)^{-\frac{3}{2}\left(\frac{1}{q}-\frac{1}{r}\right)-\frac{j}{2}}\|\rho^{\alpha}g\|_{q, \mathbb R^3}
= C(t-s)^{-\frac{3}{2}\left(\frac{1}{q}-\frac{1}{r}\right)-\frac{j}{2}}\|\rho^{\alpha}f\|_{q, \mathbb R^3}.
\end{align*}
Therefore we have the following; suppose that $\eta=0$, $\alpha\in[0, 3)$ and $(q, r)\in\mathcal{E}(\alpha)$ in additions to the assumptions of Lemma \ref{so-far-whole}. Then there is a constant $C=C(q, r, \alpha, j)>0$ such that
\begin{align}\label{wh-st-rot}
\|\rho^{\alpha}\nabla^{j}U(t, s)f\|_{r, \mathbb R^3}
\le C(t-s)^{-\frac{j}{2}-\frac{3}{2}\left(\frac{1}{q}-\frac{1}{r}\right)}\|\rho^{\alpha}f\|_{q, \mathbb R^3}
\end{align}
for all $t>s\ge 0$. 
\end{remark}

\subsection{Interior flow, $\{V(t, s)\}_{t\ge s\ge 0}$}\label{Vpre}
Let $1<q<\infty$ and $R>1$. We now consider the non-autonomous system
\begin{align*}
&\partial_tu+\mathcal L(t)u+\nabla p=0,\qquad \mathrm{div}\,u=0\qquad \mbox{in}\;D\times (s, \infty)\\
&u|_{\partial D_R}=0,\qquad u(\cdot, s)=f,
\end{align*}
where $\mathcal L(t)$ is the differential operator defined in \eqref{pre-proj-0} and $D_R:=D\cap B_R$. 
Applying the Fujita-Kato projection associated with the Helmholtz decomposition, we define the family of linear operators $\{L_{R}(t)\}_{t\ge0}$ on $L^q_{\sigma}(D)$ as
\begin{align}\label{ig-bdd}
\begin{split}
D_q(L_R(t))&=L^q_{\sigma}(D_R)\cap W^{1, q}_0(D_R)\cap W^{2, q}(D_R)\\
L_R(t)u&=P_{D_R}{\mathcal L}(t)u=-P_{D_R}\Delta u-(\eta(t)+\omega(t)\times x)\cdot\nabla u+\omega(t)\times u.
\end{split}
\end{align}
Since we have the property \eqref{grund}, $P_{D_R}$ is not needed in the lower order terms of \eqref{ig-bdd}. 
For this interior problem, one can apply the Tanabe-Sobolevski\u\i\, theory to $\{L_R(t)\}_{t\ge0}$ and find the evolution operator $\{V(t, s)\}_{t\ge s\ge 0}$ of parabolic type on $L^q_{\sigma}(D_R)$. 
Moreover, the following estimates were shown in \cite{Hi18, Hi20}, see also \cite[Lemma 3.7]{AH25}. 
\begin{lemma}[\cite{AH25, HR14, Hi18, Hi20}]
Suppose that $\eta$ and $\omega$ fulfill \eqref{ass2-rigid} for some $\theta\in (0, 1]$. Let $1<q\le r\le \infty, q\neq \infty$. For each $\mathcal T\in (0, \infty), m\in (0, \infty)$ and $j=0, 1$, there are constants $C_j=C_j(\mathcal T, m, q, r, \theta, D_R)>0$ and $C_2=C_2(\mathcal T, m, q, \theta, D_R)>0$ such that
\begin{align}
&\|\nabla^j V(t,s)f\|_{r,D_R}\leq C_j(t-s)^{-\frac{j}{2}-\frac{3}{2}\left(\frac{1}{q}-\frac{1}{r}\right)}\|f\|_{q,D_R} \label{sm-bdd-unif}\\
&\|p(t)\|_{q,D_R}+\|\partial_tV(t,s)f\|_{W^{-1,q}(D_R)}\leq C_2 (t-s)^{-\frac{1}{2}\left(1+\frac{1}{q}\right)}\|f\|_{q,D_R} \label{pressure-bdd-unif}
\end{align}
for all $(t, s)$ with $t-s\le \mathcal T$ as well as $0\le s<t$ and $f\in L^q_{\sigma}(D_R)$ whenever \eqref{ever} is satisfied, where $p(t)$ is the pressure associated with $V(t, s)f$ and it is singled out subject to the side condition $\int_{D_R}p\,dx=0$.  
\end{lemma}

\subsection{Uniform estimate in $t-s$}
For the proof of Theorem \ref{thm1}, we need the smoothing estimates of the evolution operator $T(t, s)$ uniformly in the difference $t-s$. 
Hishida and the present author \cite[Theorem 2.1]{AH25} proved the weighted smoothing estimates of $T(t, s)$ for $0\le s<t\le \mathcal T$, however, such estimates also hold for all $(t, s)$ with $t-s\le \mathcal T$, where $\mathcal T\in (0, \infty)$ is arbitrarily given, under the reasonable assumption \eqref{ass2-rigid}. 
The proof of the next proposition is the same as in \cite[Section 5]{AH25} with the aid of \eqref{wei-wh-unif}--\eqref{weak-wh-unif} and \eqref{sm-bdd-unif}--\eqref{pressure-bdd-unif}, so we omit the details.
Note that Hishida \cite[Proposition 3.1]{Hi18} proved the following proposition when $\alpha=0$.
\begin{proposition}\label{unifest-wei}
Suppose that $\eta$ and $\omega$ fulfill \eqref{ass2-rigid} for some $\theta\in (0, 1]$. Let $\alpha\in [0, 3), (q, r)\in \mathcal E(\alpha)$ and $j=0, 1$. For each $\mathcal T\in (0, \infty)$ and $m\in (0, \infty)$, there is a constant $C=C(\mathcal T, m, \alpha, q, r, \theta, D)>0$ such that
\begin{align}\label{est-uni-wei}
\|\rho^{\alpha}\nabla^jT(t, s)f\|_r\le C(t-s)^{-\frac{3}{2}\left(\frac{1}{q}-\frac{1}{r}\right)-\frac{j}{2}}\|\rho^{\alpha}f\|_{q}
\end{align}
for all $(t, s)$ with 
\begin{align*}
\mbox{$t-s\le \mathcal T$ as well as $0\le s<t$}
\end{align*}
and $f\in L^q_{\alpha, \sigma}(D)$ whenever \eqref{ever} is satisfied.
\end{proposition}
The other estimate uniformly in $t-s$ regarding the modest regularity is given by the following proposition.
For the proof, see \cite[Section 5]{Hi20}.
\begin{proposition}[\cite{Hi20}]
Let $1<q<\infty$. Given $f\in L^q_{\sigma}(D)$, we set $u(t)=T(t, s)f$. For each $\mathcal T\in (0, \infty)$ and $m\in (0, \infty)$, there is a constant $C=C(\mathcal T, m, q, \theta, D)>0$ such that
\begin{align}
&u\in C^{1}((s, \infty); W^{-1, q}(D_3)),\notag \\
&\|\partial_tT(t, s)f\|_{W^{-1, q}(D_3)}\le C(t-s)^{-\frac{1}{2}\left(1+\frac{1}{q}\right)}\|f\|_q\label{der-neg-smooth}
\end{align}
for all $(t, s)$ with $t-s\le \mathcal T$ as well as $0\le s<t$. Furthermore, we have the pressure $p(t)$ subject to $\int_{D_3}p\,dx=0$ such that the pair $\{u, p\}$ satisfies
\begin{align}\label{der-neg}
\langle \partial_tu, \psi\rangle_{D_3}+\langle\nabla u+u\otimes(\eta+\omega\times x)-(\omega\times x)\otimes u, \nabla\psi\rangle_{D_3}-\langle p, \mathrm{div}\;\psi\rangle_{D_3}=0
\end{align}
for all $t\in (s, \infty)$ and $\psi\in W^{1, q'}_{0}(D_3)$, that
\begin{align}\label{led-pre-neg}
\|p(t)\|_{q, D_3}\le C\|\partial_tu(t)\|_{W^{-1, q}(D_3)}+C\|u(t)\|_{W^{1, q}(D_3)}
\end{align}
for all $t\in (s, \infty)$ with a constant $C=C(m, q, D)>0$ and that
\begin{align*}
\|p(t)\|_{q, D_3}\le C(t-s)^{-\frac{1}{2}\left(1+\frac{1}{q}\right)}\|f\|_q
\end{align*}
for all $(t, s)$ with $t-s\le \mathcal T$ as well as $0\le s<t$ with a constant $C=C(\mathcal T, m, q, \theta, D)>0$, whenever \eqref{ever} is satisfied.

\end{proposition}

\section{Local energy decay estimates with weighted data}\label{sect-led}
In this section, we develop decay estimates of $T(t, s)f$ near the boundary, so called the local energy decay estimates. 
They play an important role in deduction of the large time behavior of the evolution operator. 
Regarding this issue, Hishida \cite[Propositions 6.1, 6.2]{Hi20} showed the following proposition.
Remark that the second assertion of the following proposition is a slight improvement of \cite[Proposition 6.1]{Hi20}, in which $f\in W^{1, q}(D)$ was additionally assumed.
This condition can be removed with the aid of the item 2 of Proposition \ref{prop1-gen-sm}.
\begin{proposition}\label{led-hi}
Suppose that $\eta$ and $\omega$ fulfill \eqref{ass2-rigid} for some $\theta\in (0, 1]$. 
\begin{enumerate}
\item Let $1<q<\infty$ and $\varepsilon>0$ be arbitrarily small. For each $m\in (0, \infty)$, there is a constant $C=C(m, \varepsilon, q, \theta, D)>0$ such that
\begin{align*}
\|T(t, s)f\|_{W^{1, q}(D_3)}\le C(t-s)^{-\frac{3}{2}+\varepsilon}\|f\|_q
\end{align*}
for all $(t, s)$ with $t-s>2$ as well as $0\le s<t$ and $f\in L^q_{\sigma}(D)$ with
\begin{align*}
f(x)=0\quad\mbox{a.e.}\;\mathbb R^3\setminus B_3
\end{align*}
whenever \eqref{ever} is satisfied. 
\item Let $3/2<q<\infty$ and $\varepsilon>0$ be arbitrarily small. For each $m\in (0, \infty)$, there is a constant $C=C(m, \varepsilon,  q, \theta, D)>0$ such that
\begin{align}
\|\partial_tT(t, s)f\|_{W^{-1, q}(D_3)}\le C(t-s)^{-\frac{3}{2}+\varepsilon}\|f\|_q
\label{led-weak-cpt}
\end{align}
for all $(t, s)$ with $t-s>2$ as well as $0\le s<t$ and $f\in L^q_{\sigma}(D)$ with
\begin{align}\label{ass-ae0}
f(x)=0\quad\mbox{a.e.}\;\mathbb R^3\setminus B_3
\end{align}
whenever \eqref{ever} is satisfied. 
\item Let $1<q<\infty$. For each $m\in (0, \infty)$, there is a constant $C=C(m, q, \theta, D)>0$ such that
\begin{align}
\|T(t, s)f\|_{W^{1, q}(D_3)}+\|\partial_tT(t, s)f\|_{W^{-1, q}(D_3)}\le C(t-s)^{-\frac{3}{2q}}\|f\|_q
\label{led-gen}
\end{align}
for all $(t, s)$ with $t-s>2$ as well as $0\le s<t$ and $f\in L^q_{\sigma}(D)$ whenever \eqref{ever} is satisfied. Here, the temporal derivative is understood as in \eqref{der-neg}
\end{enumerate}
\end{proposition}

\noindent
{\it Proof of the second assertion}.
Since \eqref{ass-ae0} implies $f\in L^q_{1, \sigma}(D)$, we have \eqref{st-cl} and \eqref{str} in $L^q_{\sigma}(D)$ for such $f$, see \cite[Theorem 2.2]{AH25}. 
Then the same argument as in \cite[Proposition 6.1]{Hi20} works well and leads to \eqref{led-weak-cpt} without assuming more regularity of $f$. 
\hfill
$\Box$
\medskip

The aim of this section is to develop further local energy decay estimates with faster rate when the initial data are taken from the weighted solenoidal space $L^q_{\alpha, \sigma}(D)$. 
This was already observed by Tomoki Takahashi \cite[Section 6]{T24}, who considered the Oseen semigroup with anisotropic weights.
Our results may be thus regarded as generalization of his to the non-autonomous case. 
To begin with, we consider better integrability of $f$ at infinity when $f\in L^q_{\alpha}(D)$. 
The following lemma is found in \cite[Lemma 5.2]{T24}.

\begin{lemma}\label{lem-low-int}
Let $\alpha\in [0, 3)$ and $3/(3-\alpha)<q<\infty$.
For $f\in L^q_{\alpha}(D)$, we have $f\in L^{\frac{3q}{3+\alpha q}, q}(D)$ and the estimate 
\begin{align}\label{lem-est-low-int}
\|f\|_{\frac{3q}{3+\alpha q}, q}\le C\|\rho^{\alpha}f\|_q
\end{align}
with some constant $C=C(q, \alpha, D)$, where the space $L^{q, r}(D)$ denotes the Lorentz space.
\end{lemma}

\begin{proof}
The Lorentz-H\"older inequality immediately implies 
\begin{align*}
\|f\|_{\frac{3q}{3+\alpha q}, q}\le C\|\rho^{-\alpha}\|_{\frac{3}{\alpha}, \infty}\|\rho^{\alpha} f\|_{q},
\end{align*}
which claims the assertion.
\end{proof}
With Lemma \ref{lem-low-int} at hand, we can prove the following lemma.
\begin{lemma}\label{lor-lplq-ext}
Suppose $\eta$ and $\omega$ fulfill \eqref{ass2-rigid} for some $\theta\in (0, 1]$. 
Assume that $\alpha\in [0, 3)$ and $(q, r)\in\mathcal{E}(\alpha)$. For each $m\in (0, \infty)$, there is a constant $C=C(m, \alpha, q, r, \theta, D)>0$ such that
\begin{align}\label{lorentz-lplq-ext}
\|T(t, s)f\|_{r}\le C(t-s)^{-\frac{3}{2}\left(\frac{1}{q}-\frac{1}{r}\right)-\frac{\alpha}{2}}\|\rho^{\alpha}f\|_{q}
\end{align}
for all $0\le s<t$ and $f\in L^{q}_{\alpha, \sigma}(D)$.
\end{lemma}

\begin{proof}
Since $3q/(3+\alpha q)>1$, we can take $(q_0, q_1)$ with $1<q_0<3q/(3+\alpha q)<q_1<\infty$ and $\kappa\in (0, 1)$ such that
\begin{align*}
L^{\frac{3q}{3+\alpha q}, q}_{\sigma}(D)=\left(L^{q_0}_{\sigma}(D), L^{q_1}_{\sigma}(D)\right)_{\kappa, q},\qquad \frac{\alpha}{3}+\frac{1}{q}=\frac{1-\kappa}{q_0}+\frac{\kappa}{q_1}.
\end{align*}
By interpolation with use of \eqref{lqlr-hi-0}, we have
\begin{align*}
\|T(t, s)\|_{\mathcal L\left(L^{\frac{3q}{3+\alpha q}, q}_{\sigma}(D),\, L^r_{\sigma}(D)\right)}
&\le \|T(t, s)\|^{1-\kappa}_{\mathcal L\left(L^{q_0}_{\sigma}(D),\, L^r_{\sigma}(D)\right)}\|T(t, s)\|^{\kappa}_{\mathcal L\left(L^{q_1}_{\sigma}(D),\, L^r_{\sigma}(D)\right)}\\
&\le C(t-s)^{-\frac{3}{2}\left(\frac{1-\kappa}{q_0}+\frac{\kappa}{q_1}-\frac{1}{r}\right)}\\
&\le C(t-s)^{-\frac{3}{2}\left(\frac{1}{q}-\frac{1}{r}\right)-\frac{\alpha}{2}}.
\end{align*}
Combining this estimate with \eqref{lem-est-low-int} concludes
\begin{align*}
\|T(t, s)f\|_{r}
&\le C(t-s)^{-\frac{3}{2}\left(\frac{1}{q}-\frac{1}{r}\right)-\frac{\alpha}{2}}\|f\|_{\frac{3q}{3+\alpha q}, q}\\
&\le C(t-s)^{-\frac{3}{2}\left(\frac{1}{q}-\frac{1}{r}\right)-\frac{\alpha}{2}}\|\rho^{\alpha}f\|_{q}.
\end{align*}
The proof is complete.
\end{proof}

Once we have \eqref{lorentz-lplq-ext}, we can show the following local energy decay estimate. 
Compared to \eqref{led-gen}, the rate of decay is $(-\alpha/2)$-faster when the initial data belong to the weighted space. 

\begin{proposition}\label{wei-led}
Let $\eta$ and $\omega$ fulfill \eqref{ass2-rigid} for some $\theta\in(0, 1]$. 
Assume that $\alpha\in [0, 3)$ and $q\in(\frac{3}{3-\alpha}, \infty)$. For each $m\in (0, \infty)$, there is a constant $C=C(m, q, \theta, \alpha)>0$ such that
\begin{align}
\|T(t, s)f\|_{W^{1, q}(D_3)}+\|\partial_tT(t, s)f\|_{W^{-1, q}(D_3)}\le C(t-s)^{-\frac{3}{2q}-\frac{\alpha}{2}}\|\rho^{\alpha}f\|_q
\label{led-wei}
\end{align}
for all $(t, s)$ with $t-s>2$ as well as $0\le s<t$ and $f\in L^q_{\alpha, \sigma}(D)$ whenever \eqref{ever} is satisfied.
\end{proposition}

\begin{proof}
Since $C^{\infty}_{0, \sigma}(D)$ is dense in $L^q_{\alpha, \sigma}(D)$ when $q\in (\frac{3}{3-\alpha}, \infty)$, it suffices to show \eqref{led-wei} for all $f\in C^{\infty}_{0, \sigma}(D)$.
Then $T(t, s)f$ is of class $C^1$ in time with values in $L^r_{\sigma}(D)$ provided $r\in (3/2, \infty)$, see Proposition \ref{prop1-gen-sm}.
Even if $q\le 3/2<r$, we still have the relation \eqref{str} in $W^{-1, r}(D_3)\subset W^{-1, q}(D_3)$.
Then we use the semigroup property to combine \eqref{led-gen} with \eqref{lorentz-lplq-ext};
as a consequence, we conclude
\begin{align*}
\|\partial_t T(t, s)f\|_{W^{-1, q}(D_3)}
&=\|L(t)T(t, (t+s)/2)T((t+s)/2, s)f\|_{W^{-1, q}(D_3)}\\
&\le C(t-s)^{-\frac{3}{2q}}\|T((t+s)/2, s)f\|_{q}\\
&\le C(t-s)^{-\frac{3}{2q}-\frac{\alpha}{2}}\|\rho^{\alpha}f\|_q.
\end{align*}
for $t-s>2$.
In the similar fashion, one can derive the decay property in $W^{1, q}(D_3)$.
The proof is complete.
\end{proof}

For later use, it is convenient to prepare the following corollary.
Estimate \eqref{led-dt-p} was already found in \cite[Corollary 6.1]{Hi20}.
Similarly, with the aid of \eqref{led-wei}, the decay estimate \eqref{led-dt-p-wei} with $\alpha\in [0, 3)$ is available. 

\begin{corollary}\label{prop-led-dt-p}
Suppose that $\eta$ and $\omega$ fulfill \eqref{ass2-rigid} for some $\theta \in (0, 1]$.
Let $\alpha\in [0, 3)$ and $q\in(\frac{3}{3-\alpha}, \infty)$.
Given $f\in L^q_{\alpha, \sigma}(D)$, we denode by $p(t)$ the pressure associated with $T(t, s)f$ subject to the condition $\int_{D_3}p\,dx=0$.
Let $\phi \in C^{\infty}_{0}(B_3)$ satisfy $\phi\equiv1$ in $B_2$ and $\mathbb B$ is the Bogovski\u\i\, operator on the annulus $A_3:=\{1<|x|<3\}$, see Lemma \ref{lem-bog}.
Then, for each $m\in (0, \infty)$, there are constants $C_0=C_0(m, q, R, \theta, D)>0$ and $C_1=C_1(m, q, R, \theta, D, \alpha)>0$ such that
\begin{align}
&\|p(t)\|_{q, D_3}+\|\mathbb B[(\partial_tT(t, s)f)\cdot\nabla \phi]\|_{q, A_3}\le C_0\|f\|_q
\begin{cases}
(t-s)^{-\frac{1}{2}\left(1+\frac{1}{q}\right)}, \quad&0<t-s\le2,\\
(t-s)^{-\frac{3}{2q}},&t-s>2
\end{cases}\label{led-dt-p}
\\
&\|p(t)\|_{q, D_3}+\|\mathbb B[(\partial_tT(t, s)f)\cdot\nabla \phi]\|_{q, A_3}\le C_1(t-s)^{-\frac{3}{2q}-\frac{\alpha}{2}}\|\rho^{\alpha}f\|_q,\qquad \;t-s>2,\label{led-dt-p-wei}
\end{align}
whenever \eqref{ever} is satisfied.
Here, the temporal derivative is understood as in \eqref{der-neg}.
\end{corollary}

\begin{proof}
Since $\nabla \phi$ vanishes near $\partial A_3$, the estimate \eqref{bog-est2} leads to
\begin{align*}
\|\mathbb B[(\partial_tT(t, s)f)\cdot\nabla \phi]\|_{q, A_3}
&\le C\|(\partial_tT(t, s)f)\cdot\nabla \phi\|_{W^{1, q'}(A_3)^{*}}\\
&\le C\|\partial_tT(t, s)f\|_{W^{-1, q}(A_3)}.
\end{align*}
From this along with \eqref{led-pre-neg} it follows that \eqref{led-wei} and \eqref{der-neg-smooth} implies \eqref{led-dt-p}--\eqref{led-dt-p-wei}.
The proof is complete. 
\end{proof}

\section{Proof of Theorem \ref{thm1}}\label{sect-dec}
To show the large time behevior of $\rho^{\alpha}\nabla^jT(t, s)$, we have to establish the weighted $L^q$-$L^r$ estimates near spatial infinity. 
Once we have such estimates, we can obtain \eqref{dec-0-alpha1}--\eqref{dec-1-alpha2} by combining those with the local energy decay estimates \eqref{led-gen} and \eqref{led-wei}. 
When $\alpha\in (0, 1)$, estimate \eqref{led-wei} plays an important role to enlarge the upper bound for $q$ to obtain the desired rate as in the right-hand side of \eqref{rmk-thm1-lqlq} in Remark \ref{rmk-thm1-1}.
\par
Let $\alpha\in [0, 1)$, then it suffices to show 
\begin{gather}
\|\rho^{\alpha}T(t, s)f\|_{r, D^c_3}
\le C(t-s)^{-\frac{3}{2}\left(\frac{1}{q}-\frac{1}{r}\right)}\|\rho^{\alpha}f\|_q+C(t-s)^{-\frac{3}{2}\left(\frac{1}{q}-\frac{1}{r}\right)+\alpha}\|f\|_q	,\label{dec-0-infty-01}\\
\|\rho^{\alpha}\nabla T(t, s)f\|_{q, D^c_3}
\le C(t-s)^{-\frac{1}{2}}\|\rho^{\alpha}f\|_q+ C(t-s)^{\max\left\{-\frac{1}{2}+\alpha, -\frac{3}{2q}-\frac{\alpha}{2}\right\}}\|f\|_{q}\label{dec-1-infty-01}
\end{gather}
for all $t-s>2$ and $f\in C^{\infty}_{0, \sigma}(D)$ since $C^{\infty}_{0, \sigma}(D)$ is dense in $L^q_{\alpha, \sigma}(D)$ when $q\in (\frac{3}{3-\alpha}, \infty)$, where $r$ satisfies $1/q-1/r<2/3$.
Given $f\in C_{0, \sigma}^{\infty}(D)$, we set $u(t):=T(t, s)f$.
We denote the pressure associated with $u(t)$ by $p(t)$, that is singled out in such a way that the condition
\begin{align*}
\int_{D_3}p(t)\,dx=0
\end{align*} 
is satisfied for all $t>0$.
We fix a cut-off function $\phi\in C^{\infty}_{0}(B_3)$ which is identically equal to 1 in $B_2$. 
We define the pair $(v, p_v)$ by 
\begin{align*}
&v(t):=(1-\phi)T(t, s)f+\mathbb B[(\nabla \phi)\cdot T(t, s)f],\\
&p_v(t):=(1-\phi)p(t),
\end{align*}
where $\mathbb B$ is the Bogovski\u\i\, operator on the annulus $A_3=\{1<|x|<3\}$.
Since $v(t)=T(t, s)f$ for $|x|\ge 3$, we are going to estimate $\rho^{\alpha}\nabla^jv(t)$ ($j=0, 1$). 
The pair of $v$ and $p$ obeys
\begin{align}\label{v-eq}
\begin{split}
&\partial_tv+\mathcal L(t)v+\nabla p_v=H,\quad \mathrm{div}\;v=0, \qquad \mbox{in}\; \mathbb R^3\times (0, \infty), \\
&v(s)=\tilde{f}:=(1-\phi)f+\mathbb B[(\nabla \phi)\cdot f]\in L^q_{\alpha, \sigma}(\mathbb R^3),
\end{split}
\end{align}
where the remainder term $H(t)$, which has a compact support, is defined by
\begin{align*}
H:=&2(\nabla \phi)\cdot\nabla T(t, s)f+\{\Delta \phi+(\eta+\omega\times x)\cdot\nabla \phi\}T(t, s)f\\
&-\Delta\mathbb B[(\nabla\phi)\cdot T(t, s)f]-(\eta+\omega\times x)\cdot\nabla \mathbb B[(\nabla\phi)\cdot T(t, s)f]\\
&+\omega\times\mathbb B[(\nabla \phi)\cdot T(t, s)f]+\mathbb B[(\nabla \phi)\cdot\partial_t T(t, s)f]-(\nabla \phi)p.
\end{align*}
Here, it follows from \eqref{led-gen} and \eqref{led-dt-p} together with \eqref{lqlr-hi-0}--\eqref{lqlr-hi-1} that
\begin{align}\label{est-H-normal}
\|H(t)\|_{\tilde{r}, \mathbb R^3}\le C(m+1)\|f\|_q
\begin{cases}
(t-s)^{-\frac{1}{2}\left(1+\frac{1}{q}\right)},\quad &0<t-s\le 2\\
(t-s)^{-\frac{3}{2q}},&t-s>2
\end{cases} 
\end{align}
for all $(t, s)$ with $0\le s< t$ and $\tilde{r}\in (1, q]$.
Moreover, on account of the local energy decay estimates \eqref{led-wei} and \eqref{led-dt-p-wei} with weighted data, estimate \eqref{est-H-normal} for $t-s>2$ can be replaced by
\begin{align}\label{est-H-wei}
\|H(t)\|_{\tilde{r}, \mathbb R^3}\le C(t-s)^{-\frac{3}{2q}-\frac{\alpha}{2}}\|\rho^{\alpha}f\|_q
\end{align}
with $\tilde{r}\in (1, q]$. 
By appling the Fujita-Kato projection $P_{\mathbb R^3}$ to the equation \eqref{v-eq}, it is converted into the integral equation
\begin{align}\label{dec-vdef}
v(t)=U(t, s)\tilde{f}+\overline{v}(t),\qquad \overline{v}(t):= \int_{s}^{t}U(t, \tau)P_{\mathbb R^3}H(\tau)\,d\tau.
\end{align}

We employ \eqref{wei-wh-all} to find
\begin{align}\label{dec-wh}
\|\rho^{\alpha}\nabla^jU(t, s)\tilde{f}\|_{r, \mathbb R^3}
\le C(t-s)^{-\frac{j}{2}-\frac{3}{2}\left(\frac{1}{q}-\frac{1}{r}\right)}\|\rho^{\alpha}f\|_{q}
+C(t-s)^{-\frac{j}{2}-\frac{3}{2}\left(\frac{1}{q}-\frac{1}{r}\right)+\alpha}\|f\|_{q}
\end{align}
for all $0\le s<t$ and $j=0, 1$.
Our task is thus to estimate $\|\rho^{\alpha}\nabla^{j}\overline{v}(t)\|_{r, \mathbb R^3}\;$for $(q, r)\in\mathcal{E}(\alpha)$ and $j=0, 1$.
Let us state an overview of the argument. 
For the case $0\le \alpha<1$, we follow almost the same calculation as in \cite[Section 7]{T24}, where we need to take an auxiliary exponent $r_0$ satisfying \eqref{r0cond} below.
Here we emphasize that we should use \eqref{est-H-normal} rather than \eqref{est-H-wei} to obtain the $L^q$-$L^r$ estimates of the desired form \eqref{dec-0-alpha1}--\eqref{dec-1-alpha1}.\par
As for the case $1\le\alpha<3$, this procedure cannot work because the auxiliary exponent $r_0$ mentioned above does not exist.  
To overcome this difficulty, we give up deducing the estimates of the form \eqref{dec-0-alpha1}--\eqref{dec-1-alpha1}, but instead, our idea is to apply \eqref{est-H-wei} toward \eqref{dec-0-alpha2}--\eqref{dec-1-alpha2}.
More precisely, we choose another auxiliary exponent satisfying \eqref{r1cond} for the case $1\le\alpha<2$, whereas we do not have to take any auxiliary exponent otherwise (case $2\le\alpha<3$).
This procedure allows us to extend the range of $\alpha$ for obtaining alternative estimates \eqref{dec-0-alpha2}--\eqref{dec-1-alpha2}. 
Compared with the case $0\le\alpha<1$, however, the large time behavior becomes worse when $1\le\alpha<3$. \par
Now, we split $\overline{v}(t)$ into three parts as follows:
\begin{align*}
\overline{v}(t)=\left(\int_{s}^{\frac{s+t}{2}}+\int_{\frac{s+t}{2}}^{t-1}+\int_{t-1}^{t}\right)U(t, \tau)P_{\mathbb R^3}H(\tau)\,d\tau=:\overline{v_1}(t)+\overline{v_2}(t)+\overline{v_3}(t)
\end{align*}

The following lemma provides the key estimates.

\begin{lemma}\label{lem-flow-infty}
Suppose that $\eta$ and $\omega$ fulfill \eqref{ass2-rigid} for some $\theta\in (0, 1]$.
Assume that $\alpha\in [0, 3)$ and $(q, r)\in\mathcal{E}(\alpha)$. 
\begin{enumerate}
\item Let $0\le \alpha<1$. We further assume that $1/q-1/r<2/3$ for the estimate \eqref{vb-0-01} below.
For each $m\in (0, \infty)$, there are constants $C_1=C_1(m, q, r, \theta, \alpha, D)>0$ and $C_2=C_2(m, q, \theta, \alpha, D)>0$ such that
\begin{align}
&\|\rho^{\alpha}\overline{v}(t)\|_r\le C_1(t-s)^{-\frac{3}{2}\left(\frac{1}{q}-\frac{1}{r}\right)+\alpha}\|f\|_q\label{vb-0-01}\\
&\|\rho^{\alpha}\nabla\overline{v_1}(t)\|_{q, \mathbb R^3}\le C_2(t-s)^{-\frac{1}{2}+\alpha}\|f\|_q\label{vb1-1-01}\\
&\|\rho^{\alpha}\nabla\overline{v_2}(t)\|_{q, \mathbb R^3}\le C_2(t-s)^{\max\left\{-\frac{1}{2}+\alpha, -\frac{3}{2q}\right\}}\|f\|_q\label{vb2-1-01}\\
&\|\rho^{\alpha}\nabla\overline{v_3}(t)\|_{q, \mathbb R^3}\le C_2(t-s)^{-\frac{3}{2q}}\|f\|_q\label{vb3-1-01}
\end{align}
for all $(t, s)$ with $t-s>2$  as well as $0\le s<t$, whenever \eqref{ever} is satisfied. 

\item Let $1\le\alpha<3$ and $j=0, 1$. 
We further assume $1/q-1/r<1/3$ for \eqref{vb3-j-13} below when $1\le\alpha< 2$.
For each $m\in (0, \infty)$, there are constants $C^{\prime}_k=C^{\prime}_{k}(m, q, r, \theta, \alpha, D)>0\, (k=1, 2)$  such that
\begin{align}
\|\rho^{\alpha}\nabla^j\overline{v_1}(t)\|_{r, \mathbb R^3}+\|\rho^{\alpha}\nabla^j\overline{v_2}(t)\|_{r, \mathbb R^3}
&\le C^{\prime}_1(t-s)^{-\frac{j}{2}-\frac{3}{2}\left(\frac{1}{q}-\frac{1}{r}\right)+\alpha}\|\rho^{\alpha}f\|_q\label{vb12-j-13}\\
\|\rho^{\alpha}\nabla^j\overline{v_3}(t)\|_{r, \mathbb R^3}
&\le C^{\prime}_2
(t-s)^{-\frac{3}{2q}}\|f\|_q\label{vb3-j-13}
\end{align} 
for all $(t, s)$ with $t-s>2$ as well as $0\le s<t$, whenever \eqref{ever} is satisfied.
\end{enumerate}
\end{lemma}

\begin{proof}
Throughout the proof, every estimate holds for all $(t, s)$ with $t-s>2$ as well as $0\le s<t$, although this will not be always mentioned explicitly.
\par
By using \eqref{wei-wh-all} together with \eqref{est-H-normal}--\eqref{est-H-wei}, and by taking into account
\begin{align*}
\|\rho^{\alpha}P_{\mathbb R^3}H(\tau)\|_{\tilde{r}, \mathbb R^3}\le C\|H(\tau)\|_{\tilde{r}, \mathbb R^3}
\end{align*}
for $3/(3-\alpha)<\tilde r$, see \eqref{FK-wei}, we have 
\begin{align}\label{def-kj}
\|\rho^{\alpha}\nabla^jU(t, \tau)P_{\mathbb R^3}H(\tau)\|_{r, \mathbb R^3}
\le C
\begin{cases}
K_j(\tau)\|f\|_q,\\
K'_j(\tau)\|\rho^{\alpha}f\|_q,
\end{cases}
\end{align}
where
\begin{align}
K_j(\tau):=&(t-\tau)^{-\frac{j}{2}-\frac{3}{2}\left(\frac{1}{q}-\frac{1}{r}\right)}(1+t-\tau)^{-\frac{3}{2}\left(\frac{1}{\tilde{r}}-\frac{1}{q}\right)+\alpha}(\tau-s)^{-\frac{1}{2}\left(1+\frac{1}{q}\right)}(1+\tau-s)^{\frac{1}{2}-\frac{1}{q}},\label{dec-0-01}\\
K'_j(\tau):=&(t-\tau)^{-\frac{j}{2}-\frac{3}{2}\left(\frac{1}{q}-\frac{1}{r}\right)}(1+t-\tau)^{-\frac{3}{2}\left(\frac{1}{\tilde{r}}-\frac{1}{q}\right)+\alpha}(\tau-s)^{-\frac{1}{2}\left(1+\frac{1}{q}\right)}
(1+\tau-s)^{\frac{1}{2}-\frac{1}{q}-\frac{\alpha}{2}},\label{dec-1-01}
\end{align}
for $\alpha\in [0, 3), (q, r)\in \mathcal{E}(\alpha)$ and $\tilde{r}\in (\frac{3}{3-\alpha}, q]$.
\par
We first verify \eqref{vb-0-01}--\eqref{vb3-1-01}.
On account of $0\le \alpha<1$, we can take an exponent $r_0$ satisfying
\begin{align}\label{r0cond}
\frac{3}{3-\alpha}<r_0<\min\left\{\frac{3}{2}, q\right\}.
\end{align}
Let us begin with estimate of $\rho^{\alpha}\nabla^j\overline{v_1}$.
By \eqref{def-kj}--$\eqref{dec-0-01}_{j=0}$ with $\tilde r=r_0$, 
\begin{align*}
\|\rho^{\alpha}\overline{v_1}(t)\|_{r, \mathbb R^3}
&\le C(t-s)^{-\frac{3}{2}\left(\frac{1}{r_0}-\frac{1}{r}\right)+\alpha}\int_{s}^{\frac{s+t}{2}}(\tau-s)^{-\frac{1}{2}\left(1+\frac{1}{q}\right)}(1+\tau-s)^{\frac{1}{2}-\frac{1}{q}}\,d\tau\|f\|_q\\
&\le C(t-s)^{-\frac{3}{2}\left(\frac{1}{r_0}-\frac{1}{r}\right)+\alpha}R(t-s)\|f\|_q,
\end{align*}
for $t-s>2$ as well as $0\le s<t$, where
\begin{align*}
R(t-s):=
\begin{cases}
(t-s)^{1-\frac{3}{2q}},\qquad &\frac{3}{2}<q<\infty,\\
\log(t-s),\qquad & q=\frac{3}{2},\\
1,\qquad &\frac{3}{3-\alpha}<q<\frac{3}{2}.
\end{cases}
\end{align*}
It follows from \eqref{r0cond} that 
\begin{align}\label{dec-vb1-0-01}
\|\rho^{\alpha}\overline{v_1}(t)\|_{r, \mathbb R^3}\le C(t-s)^{-\frac{3}{2}\left(\frac{1}{q}-\frac{1}{r}\right)+\alpha}\|f\|_q
\end{align}
for $t-s>2$.
Similarly, \eqref{def-kj}--$\eqref{dec-0-01}_{j=1}$ with $\tilde r=r_0$ leads to \eqref{vb1-1-01}; in fact,
\begin{align*}
\|\rho^{\alpha}\nabla\overline{v_1}(t)\|_{q, \mathbb R^3}
&\le C(t-s)^{-\frac{1}{2}-\frac{3}{2}\left(\frac{1}{r_0}-\frac{1}{q}\right)+\alpha}R(t-s)\|f\|_{q}\\
&\le C(t-s)^{-\frac{1}{2}+\alpha}\|f\|_q.
\end{align*}
\par
As for $\rho^{\alpha}\overline{v_2}$, we employ \eqref{def-kj}--$\eqref{dec-0-01}_{j=0}$ with $\tilde r=r_0$ to find
\begin{align*}
\|\rho^{\alpha}\overline{v_2}(t)\|_{r, \mathbb R^3}
&\le C\|f\|_q (t-s)^{-\frac{3}{2q}}\int_{1}^{\frac{t-s}{2}}\xi^{-\frac{3}{2}\left(\frac{1}{q}-\frac{1}{r}\right)}(1+\xi)^{-\frac{3}{2}\left(\frac{1}{r_0}-\frac{1}{q}\right)+\alpha}\,d\xi,\\
&\le C(t-s)^{-\frac{3}{2}\left(\frac{1}{q}-\frac{1}{r}\right)+\alpha}\|f\|_q\int_{1}^{\frac{t-s}{2}}\xi^{-\frac{3}{2r_0}}\,d\xi,
\end{align*}
yielding
\begin{align}\label{dec-vb2-0-01}
\|\rho^{\alpha}v_2(t)\|_{r, \mathbb R^3}\le C(t-s)^{-\frac{3}{2}\left(\frac{1}{q}-\frac{1}{r}\right)+\alpha}\|f\|_q.
\end{align}
Concerning $\|\rho^{\alpha}\nabla \overline{v_2}\|_{q, \mathbb R^3}$, besides the condition \eqref{r0cond}, one can choose particularly the exponent $r_0$ such that
\begin{align*}
\frac{1}{r_0}-\frac{1}{q}\neq \frac{1+2\alpha}{3}.
\end{align*}
Indeed, we define the exponents
\begin{align*}
p_0:=\frac{q_0+q_1}{2},\qquad p_1:=\frac{q_0+2q_1}{3},
\end{align*}
to take $r_0$ by
\begin{align*}
r_0=
\begin{cases}
p_0,\qquad \frac{1}{p_0}-\frac{1}{q}\neq \frac{1+2\alpha}{3},\\
p_1,\qquad \frac{1}{p_0}-\frac{1}{q}= \frac{1+2\alpha}{3},\\
\end{cases}
\end{align*}
where $q_0:=3/(3-\alpha), q_1:=\min\{3/2, q\}$. 
Then the same calculation as before leads to
\begin{align*}
\|\rho^{\alpha}\nabla \overline{v_2}(t)\|_{q, \mathbb R^3}
&\le C(t-s)^{-\frac{3}{2q}+\max\left\{\frac{1}{2}-\frac{3}{2}\left(\frac{1}{r_0}-\frac{1}{q}\right)+\alpha, 0\right\}}\|f\|_q.\\
&\le C(t-s)^{\max\left\{-\frac{1}{2}+\alpha, -\frac{3}{2q}\right\}}\|f\|_q
\end{align*}
for $t-s>2$, therefore, we get \eqref{vb2-1-01}.\par
Let us proceed to the estimate of $\rho^{\alpha}\nabla^j\overline{v_3}$.
We employ \eqref{def-kj}--$\eqref{dec-0-01}_{j=0}$ with $\tilde r=q$ to obtain
\begin{align*}
\|\rho^{\alpha}\overline{v_3}(t)\|_{r, \mathbb R^3}
\le C(t-s)^{-\frac{3}{2q}}\int_{t-1}^{t}(t-\tau)^{-\frac{3}{2}\left(\frac{1}{q}-\frac{1}{r}\right)}(1+t-\tau)^{\alpha}\,d\tau\|f\|_q.
\end{align*}
Since $1/q-1/r<2/3$, the integral in the right-hand side converges. Therefore we get
\begin{align}\label{dec-vb3-0-01}
\|\rho^{\alpha}\overline{v_3}(t)\|_{r, \mathbb R^3}
\le C(t-s)^{-\frac{3}{2q}}\|f\|_q.
\end{align}
Collecting \eqref{dec-vb1-0-01}--\eqref{dec-vb3-0-01} leads to \eqref{vb-0-01}.
Similarly, \eqref{def-kj}--$\eqref{dec-0-01}_{j=1}$ with $r=\tilde r=q$ implies \eqref{vb3-1-01}.
\par
Now, we turn to the case $1\le \alpha<3$.
As already mentioned, we rely on \eqref{est-H-wei} rather than \eqref{est-H-normal} for $t-s>2$.
If in particular $1\le\alpha<2$, the other auxiliary exponent $r_1$ is taken in such a way that
\begin{align}\label{r1cond}
\frac{3}{3-\alpha}<r_1<\min\left\{\frac{3}{2-\alpha}, q\right\}.
\end{align}
Remark that we do not need to take such $r_1$ when $2\le \alpha<3$, because the additional decay rate $-\alpha/2$ in \eqref{dec-1-01} is so fast to get the estimate \eqref{dec-0-alpha2}--\eqref{dec-1-alpha2}.\par
When $1\le \alpha<2$, it follows from \eqref{def-kj} with \eqref{dec-1-01} for $\tilde{r}=r_1$ that
\begin{align*}
\|\rho^{\alpha}\nabla^j\overline{v_1}(t)\|_{r, \mathbb R^3}
&\le C(t-s)^{-\frac{j}{2}-\frac{3}{2}\left(\frac{1}{r_1}-\frac{1}{r}\right)+\alpha}\int_{s}^{\frac{s+t}{2}}(\tau-s)^{-\frac{1}{2}\left(1+\frac{1}{q}\right)}(1+\tau-s)^{\frac{1}{2}-\frac{1}{q}-\frac{\alpha}{2}}\,d\tau\|\rho^{\alpha}f\|_{q}\\
&\le C(t-s)^{-\frac{j}{2}-\frac{3}{2}\left(\frac{1}{r_1}-\frac{1}{r}\right)+\alpha}\,\overline{R}(t-s)\|\rho^{\alpha}f\|_q
\end{align*}
for $t-s>2$ as well as $0\le s<t$, where
\begin{align*}
\overline{R}(t-s):=
\begin{cases}
(t-s)^{1-\frac{3}{2q}-\frac{\alpha}{2}},\qquad&\frac{3}{2-\alpha}<q<\infty,\\
\log(t-s),\qquad&q=\frac{3}{2-\alpha},\\
1,\qquad&\frac{3}{3-\alpha}<q<\frac{3}{2-\alpha}.
\end{cases}
\end{align*}
Taking into account \eqref{r1cond}, we obtain
\begin{align}\label{vb-j-12}
\|\rho^{\alpha}\nabla^j\overline{v_1}(t)\|_{r, \mathbb R^3}
\le C(t-s)^{-\frac{j}{2}-\left(\frac{1}{q}-\frac{1}{r}\right)+\alpha}\|\rho^{\alpha}f\|_q
\end{align}
for $t-s>2$ provided $1\le \alpha<2$.
For the case $2\le \alpha<3$, \eqref{def-kj} with \eqref{dec-1-01} for $\tilde r=q$ leads us to
\begin{align*}
\|\rho^{\alpha}\nabla^j\overline{v_1}(t)\|_{r, \mathbb R^3}
&\le C(t-s)^{-\frac{j}{2}-\frac{3}{2}\left(\frac{1}{q}-\frac{1}{r}\right)+\alpha}\int_{s}^{\frac{s+t}{2}}(\tau-s)^{-\frac{1}{2}\left(1+\frac{1}{q}\right)}(1+\tau-s)^{\frac{1}{2}-\frac{1}{q}-\frac{\alpha}{2}}\,d\tau\|\rho^{\alpha}f\|_{q}.
\end{align*}
The integral in the right-hand side can be estimated from above by
\begin{align*}
\int_{s}^{s+1}(\tau-s)^{-\frac{1}{2}\left(1+\frac{1}{q}\right)}\,d\tau+\int_{s+1}^{\frac{s+t}{2}}(\tau-s)^{-\frac{3}{2q}-\frac{\alpha}{2}}\,d\tau\le C
\end{align*}
for all $q\in (\frac{3}{3-\alpha}, r]$, because the power $1-3/2q-\alpha/2$ is negative whenever $2\le \alpha<3$.
Therefore \eqref{vb-j-12} holds true for $2\le \alpha<3$ as well.\par
Note that $\alpha\ge1$ implies that $3/(3-\alpha)\ge3/2$ and $1/\tilde r-1/r<2/3$, yielding
\begin{align*}
1-\frac{j}{2}-\frac{3}{2}\left(\frac{1}{\tilde r}-\frac{1}{r}\right)+\alpha>1-\frac{j}{2}-\frac{3}{2}\cdot\frac{2}{3}+\alpha=-\frac{j}{2}+\alpha>0.
\end{align*}
This observation and \eqref{def-kj} with \eqref{dec-1-01} lead us to
\begin{align*}
\|\rho^{\alpha}\nabla^j\overline{v_2}(t)\|_{r, \mathbb R^3}
&\le C(t-s)^{-\frac{3}{2q}-\frac{\alpha}{2}}\int_{\frac{s+t}{2}}^{t-1}(t-\tau)^{-\frac{j}{2}-\frac{3}{2}\left(\frac{1}{\tilde r}-\frac{1}{r}\right)+\alpha}\,d\tau\|\rho^{\alpha}f\|_q\\
&\le C(t-s)^{-\frac{j}{2}-\frac{3}{2}\left(\frac{1}{q}-\frac{1}{r}\right)+\alpha}(t-s)^{1-\frac{3}{2\tilde r}-\frac{\alpha}{2}}\|\rho^{\alpha}f\|_q,
\end{align*}
where one should take $\tilde r=r_1$ for $1\le \alpha<2$, while $\tilde r=q$ for $2\le \alpha<3$.
Then we get
\begin{align*}
\|\rho^{\alpha}\nabla^j\overline{v_2}(t)\|_{r, \mathbb R^3}
\le C(t-s)^{-\frac{j}{2}-\frac{3}{2}\left(\frac{1}{q}-\frac{1}{r}\right)+\alpha}\|\rho^{\alpha}f\|_q,
\end{align*}
for $t-s>2$, which completes the proof of \eqref{vb12-j-13}.
\par
Finally, to verify \eqref{vb3-j-13}, it is enough to use \eqref{def-kj}--\eqref{dec-0-01}; in fact, we have
\begin{align*}
\|\rho^{\alpha}\nabla^j\overline{v_3}(t)\|_{r, \mathbb R^3}
\le C(t-s)^{-\frac{3}{2q}}\int_{t-1}^{t}(t-\tau)^{-\frac{j}{2}-\frac{3}{2}\left(\frac{1}{q}-\frac{1}{r}\right)}(1+t-\tau)^{\alpha}\,d\tau\|f\|_q.
\end{align*}
When $1\le\alpha<2$ and $j=0$, the integral in the right-hand side always converges for all $(q, r)\in\mathcal{E}(\alpha)$.
That is also the case when $2\le\alpha<3$ and $j=0, 1$.
As for the remaining case  $1\le \alpha<2$ and $j=1$, this integral converges provided $1/q-1/r<1/3$.
As a consequence, we have
\begin{align*}
\|\rho^{\alpha}\nabla^j\overline{v_3}(t)\|_{r, \mathbb R^3}
\le C(t-s)^{-\frac{3}{2q}}\|f\|_q.
\end{align*}
The proof is complete.
\end{proof}

Let us close this section with completion of the proof of the $L^q$-$L^r$ estimates globally in time.
\medskip

\noindent
{\it Proof of Theorem \ref{thm1}}.
By Proposition \ref{unifest-wei}, we already have the smoothing estimate uniformly in $t-s$, for instance, $t-s\le 2$.
Consider the large time behavior for $t-s>2$. 
Due to the local energy decay estimate \eqref{led-wei} together with \eqref{lqlr-hi-0}, we have
\begin{align}
&\|\rho^{\alpha}T(t, s)f\|_{r, D_3}
\le C\|T(t, s)f\|_{\infty}
\le C(t-s)^{-\frac{3}{2q}}\|f\|_q,\label{dec-prf-0}\\
&\|\rho^{\alpha}\nabla T(t, s)f\|_{q, D_3}
\le C\|T(t, s)f\|_{W^{1, q}(D_3)}
\le C(t-s)^{-\frac{3}{2q}-\frac{\alpha}{2}}\|\rho^{\alpha}f\|_q.\label{dec-prf-1}
\end{align}
for $\alpha\in [0, 3)$ and $(q, r)\in \mathcal{E}(\alpha)$.
\par
Thus, we have only to deduce \eqref{dec-0-infty-01}--\eqref{dec-1-infty-01} near spatial infinity by using $v(t)=T(t, s)f$ for $|x|\ge 3$.
We first consider the case $0\le \alpha<1$. 
By \eqref{dec-wh} and the first asserion of Lemma \ref{lem-flow-infty}, we have \eqref{dec-0-infty-01}--\eqref{dec-1-infty-01} for $t-s>2$, where $1/q-1/r<2/3$ is needed for \eqref{dec-0-infty-01}, see also the item 2 of Remark \ref{rmk-thm1-1}. 
We collect \eqref{est-uni-wei}, \eqref{dec-0-infty-01}--\eqref{dec-1-infty-01} and \eqref{dec-prf-0}--\eqref{dec-prf-1} to conclude \eqref{dec-0-alpha1} for $(q, r)\in\mathcal{E}(\alpha)$ with $1/q-1/r<2/3$, as well as \eqref{dec-1-alpha1} with $r=q\in (\frac{3}{3-\alpha}, \infty)$.
To get \eqref{dec-0-alpha1} for all $(q, r)\in\mathcal{E}(\alpha)$, we find an exponent $q_0\in (q, r)$ fulfilling $1/q-1/q_0<2/3$ and $1/q_0-1/r<2/3$. Then the semigroup property leads to
\begin{align*}
&\|\rho^{\alpha}T(t, s)f\|_r\\
\le &C(t-s)^{-\frac{3}{2}\left(\frac{1}{q_0}-\frac{1}{r}\right)}\|\rho^{\alpha}T((t+s)/2, s)f\|_{q_0}+C(t-s)^{-\frac{3}{2}\left(\frac{1}{q_0}-\frac{1}{r}\right)}(1+(t-s)/2)^{\alpha}\|T((t+s)/2, s)f\|_{q_0}\\
\le &C(t-s)^{-\frac{3}{2}\left(\frac{1}{q}-\frac{1}{r}\right)}\|\rho^{\alpha}f\|_q+C(t-s)^{-\frac{3}{2}\left(\frac{1}{q}-\frac{1}{r}\right)}(1+t-s)^{\alpha}\|f\|_q,
\end{align*}
for all $t>s\ge 0$, which completes the proof of \eqref{dec-0-alpha1}.
Similarly, \eqref{dec-1-alpha1} for $(q, r)\in\mathcal{E}(\alpha)$ follows from \eqref{dec-0-alpha1} and \eqref{dec-1-alpha1} with $r=q\in (\frac{3}{3-\alpha}, \infty)$ on account of the semigroup property.
\par
We next consider the case $1\le\alpha<3$.
From \eqref{dec-wh} and the second assertion of Lemma \ref{lem-flow-infty}, we find
\begin{gather}
\|\rho^{\alpha}\nabla^{j}T(t, s)f\|_{r, D^c_3}
\le \|\rho^{\alpha}\nabla^jv(t)\|_{r, \mathbb R^3}
\le C(t-s)^{-\frac{j}{2}-\frac{3}{2}\left(\frac{1}{q}-\frac{1}{r}\right)+\alpha}\|\rho^{\alpha}f\|_q\label{dec-j-infty-13}
\end{gather}
for $t-s>2$. 
By $\eqref{dec-j-infty-13}_{j=0}$ and \eqref{dec-prf-0} together with \eqref{est-uni-wei}, we are led to \eqref{dec-0-alpha2} for all $t>s\ge 0$.
As for the case $j=1$, it follows from \eqref{lqlr-hi-1} that
\begin{align}\label{13-local}
\|\rho^{\alpha}\nabla T(t, s)f\|_{r, D_3}
\le C\|\nabla T(t, s)f\|_r\le C(t-s)^{-\frac{1}{2}-\frac{3}{2}\left(\frac{1}{q}-\frac{1}{r}\right)+\max\left\{\frac{3}{2}\left(\frac{1}{3}-\frac{1}{r}\right), 0\right\}}\|f\|_q
\end{align}
for $t-s>2$, from which combined with $\eqref{dec-j-infty-13}_{j=1}$ and \eqref{est-uni-wei}, we conclude \eqref{dec-1-alpha2} since the rate \eqref{13-local} is enough to deduce the desired estimate on account of $\alpha\ge1$.
The proof of Theorem \ref{thm1} is complete.
\hfill
$\Box$

\begin{remark}\label{lqlr-rate-imp}
In order to obtain \eqref{dec-1-alpha1}, it is enough to deduce the $L^q$-$L^q$ estimates \eqref{vb1-1-01}--\eqref{vb3-1-01} of $\rho^{\alpha}\nabla \overline{v}(t)$.
Nevertheless, we also get the following $L^q$-$L^r$ estimates of $\rho^{\alpha}\nabla\overline v(t)$ in the same way:
\begin{align*}
\|\rho^{\alpha}\nabla \overline{v}(t)\|_{r, \mathbb R^3}\le C(t-s)^{\max\left\{-\frac{1}{2}-\frac{3}{2}\left(\frac{1}{q}-\frac{1}{r}\right)+\alpha, -\frac{3}{2q}\right\}}\|f\|_q
\end{align*}
for all $t-s>2$ provided $1/q-1/r<1/3$.
This together with the local energy decay estimate 
\begin{align*}
\|\nabla T(t, s)f\|_{r, D_3}\le C(t-s)^{-\frac{3}{2r}}\|T((s+t)/2, s)f\|_{r}\le C(t-s)^{-\frac{3}{2q}-\frac{\alpha}{2}}\|\rho^{\alpha}f\|_q,
\end{align*}
that follows from \eqref{led-gen} and \eqref{lorentz-lplq-ext},
implies 
\begin{align*}
\|\rho^{\alpha}\nabla T(t, s)f\|_r\le C(t-s)^{\max\left\{-\frac{1}{2}-\frac{3}{2}\left(\frac{1}{q}-\frac{1}{r}\right), -\frac{3}{2q}-\frac{\alpha}{2}\right\}}\|\rho^{\alpha}f\|_q+C(t-s)^{\max\left\{-\frac{1}{2}-\frac{3}{2}\left(\frac{1}{q}-\frac{1}{r}\right)+\alpha, -\frac{3}{2q}\right\}}\|f\|_q
\end{align*}
for all $t-s>2$ and $(q, r)\in \mathcal E(\alpha)$ provided $1/q-1/r<1/3$, 
the latter term of which improves the one of \eqref{dec-1-alpha1}.
However, this estimate would be less useful as we have seen in the item 2 of Remark \ref{rmk-thm1-1}.
\end{remark}

\begin{remark}\label{rem-ext-st-rot}
Let us consider the case $\eta\equiv0$, in which the linearized system is either the Stokes system with $\omega\equiv 0$ or the one around the purely rotating body.
Due to Remark \ref{rem-wh-st-rot}, the evolution operator in the whole space have the $L^q$-$L^r$ estimate \eqref{wh-st-rot} without growing rate $(t-s)^{\alpha}$. 
A similar argument to what we have done in Lemma \ref{lem-flow-infty} leads us to the following estimates: if $\alpha\in [0, 3)$ and $(q, r)\in \mathcal{E}(\alpha)$, we have
\begin{align*}
&\|\rho^{\alpha}v(t)\|_{r, \mathbb R^3}\le C(t-s)^{-\frac{3}{2}\left(\frac{1}{q}-\frac{1}{r}\right)}\|\rho^{\alpha}f\|_{q},\qquad 0\le\alpha<3,\\
&\|\rho^{\alpha}\nabla v(t)\|_{q, \mathbb R^3},\le C,
\begin{cases}
(t-s)^{\max\left\{-\frac{1}{2}, -\frac{3}{2q}-\frac{\alpha}{2}\right\}}\|\rho^{\alpha}f\|_{q},\quad 0\le\alpha<1,\\
(t-s)^{-\frac{1}{2}}\|\rho^{\alpha}f\|_q,\quad 1\le\alpha<3,
\end{cases}
\end{align*}
for $f\in L^q_{\alpha, \sigma}(D)$ and all $(t, s)$ with $t-s>2$, where $v$ is given by \eqref{dec-vdef}.
Note that $1/2\le 3/2q+\alpha/2$ provided $1\le\alpha$, hence the decay rate of $\rho^{\alpha}\nabla v(t)$ for the case $1\le \alpha<3$ is always $(t-s)^{-\frac{1}{2}}$.
From these estimates, we are led to \eqref{rot-0-imp}--\eqref{rot-1-imp}.
\end{remark}

\section{Navier-Stokes flow}\label{sect-stability}

Let us apply our linear theory to the Navier-Stokes initial value problem.
To find a unique solution to \eqref{ns-inteq}, we follow the argument \cite[Section 8]{T24}, who considered the Navier-Stokes flow in anisotropically weighted $L^q$-spaces.
\par
Before going to the proof, we give the summability of the time-periodic flow $v$ under consideration.
The following lemma is found in \cite[Lemma 5.2]{T26}.
\begin{lemma}\label{lem-stab-sum}
Let $v$ be the time-periodic solution to the \eqref{ns-per-int} obtained by Proposition \ref{prop-galtak}. Let $\lambda>0$ in \eqref{def-lambda}. Then the following statements hold:
\begin{enumerate}
\item For all $q\in (2, \infty]$, we have $v\in L^{\infty}(0, \mathcal T; L^q(D))$ with
\begin{align}\label{stab-lem-0}
\sup_{t\in(0, \mathcal T)}\|v(t)\|_q\le C\lambda^{-\frac{1}{q}}\|(\eta, \omega)\|_{W^{2, 2}(0, \mathcal T)},
\end{align}
where $C=C(q)>0$ is independent of $\lambda$.
\item For all $q\in \left(4/3, \infty\right]$, we have $\nabla v\in L^{\infty}(0, \mathcal T; L^q(D))$ with
\begin{align}\label{stab-lem-1}
\sup_{t\in(0, \mathcal T)}\|\nabla v(t)\|_q\le C\lambda^{-\frac{1}{q}}\|(\eta, \omega)\|_{W^{2, 2}(0, \mathcal T)},
\end{align}
where $C=C(q)>0$ is independent of $\lambda$.
\end{enumerate}
\end{lemma}

\begin{remark}\label{rmk-frac34}
Under the smallness assumption $\|(\eta, \omega)\|_{W^{2, 2}(0, \mathcal T)}$, one may assume that $\lambda\in (0, 1]$ and, therefore, the right-hand side of \eqref{stab-lem-0} and \eqref{stab-lem-1} is also bounded by $C\lambda^{-\frac{3}{4}}\|(\eta, \omega)\|_{W^{2, 2}(0, \mathcal T)}$.
\end{remark}

We consider the following integral equation
\begin{align}\label{inteq}
w(t)=T(t, 0)w_0-H[w, w](t)-G[w](t),
\end{align}
where
\begin{align*}
H[w_1, w_2](t):=\int_{0}^{t}T(t, \tau)P(w_1\cdot\nabla w_2)\,d\tau,\\
G[w](t):=\int_{0}^{t}T(t, \tau)P(v\cdot\nabla w+w\cdot\nabla v)\,d\tau.
\end{align*}
Let us introduce some notations. For $0\le \alpha<1$, we fix two numbers $\delta_i (i=1, 2)$ such that
\begin{align}\label{def-delta}
0<\delta_1<\min\left\{\frac{1}{2}, 1-\alpha\right\},\qquad 0<\delta_2<\min\left\{\frac{1}{4}, 1-\alpha\right\}
\end{align}
and assume that
\begin{align*}
v\in L^{\infty}\left(0, \mathcal T; L^{\frac{3}{1+\delta_1}}(D)\cap L^{\infty}(D)\right),\qquad \nabla v\in L^{\infty}\left(0, \mathcal T; L^{\frac{3}{2+\delta_2}}\cap L^3(D)\right),
\end{align*}
which is justified by Lemma \ref{lem-stab-sum}.
Here, we do need $\alpha<1$ in order that both $\delta_1$ and $\delta_2$ are less than $1-\alpha$ in \eqref{def-delta}, which stems from \eqref{FK-wei} with $f=v\cdot\nabla w+w\cdot\nabla v$.
\par
In view of the estimates \eqref{dec-0-alpha1} and \eqref{dec-1-alpha1}, we set the following norms: For $\sigma\in \{0, \alpha\}$ and $t>0$, we define 
\begin{gather*}
[w]_{q, \sigma, t}:=\sup_{0<\tau<t}\tau^{\frac{1}{2}-\frac{3}{2q}}(1+\tau)^{-\sigma}\|\rho^{\sigma}w(\tau)\|_q,\qquad [\nabla w]_{3, \sigma, t}:=\sup_{0<\tau<t}\tau^{\frac{1}{2}}(1+\tau)^{-\sigma}\|\rho^{\sigma}\nabla w(\tau)\|_3,\\
\left| \! \left| \! \left| w\right| \! \right| \! \right|_{\sigma, t}:=[w]_{3, \sigma, t}+[w]_{\infty, \sigma, t}+[\nabla w]_{3, \sigma, t},\\
\llbracket v\rrbracket:=\sup_{t\in\mathbb R}\left(\|v(t)\|_{\frac{3}{1+\delta_1}}+\|v(t)\|_{\infty}+\|\nabla v(t)\|_{\frac{3}{2+\delta_2}}+\|\nabla v(t)\|_{3}\right).
\end{gather*}

\begin{lemma}\label{stab-lem-duha}
Let $0\le\alpha<1$ and $\sigma\in\{0, \alpha\}$. 
Then there exist a constant $\mathcal{C}_{\sigma}, \widetilde{\mathcal{C}}_{\sigma}>0$ such that
\begin{align}
\left| \! \left| \! \left| H[w_1, w_2]\right| \! \right| \! \right|_{\sigma, t}&\le \mathcal{C}_{\sigma}\left([w_1]^{\frac{1}{2}}_{3, 0, t}[w_1]^{\frac{1}{2}}_{\infty, 0, t}+[w_1]^{\frac{1}{2}}_{3, \sigma, t}[w_1]^{\frac{1}{2}}_{\infty, \sigma,  t}\right)[\nabla w_2]_{3, 0, t},\label{hest}\\
\left| \! \left| \! \left| G[w]\right| \! \right| \! \right|_{\sigma, t}&\le \widetilde{\mathcal{C}}_{\sigma}\llbracket v\rrbracket \bigl([\nabla w]_{3, 0, t}+[\nabla w]_{3, \sigma, t}+[w]_{\infty, 0, t}+[w]_{\infty, \sigma, t}\bigr),\label{gest}
\end{align} 
for all $t>0$ and $w, w_j$ with $\left| \! \left| \! \left| w\right| \! \right| \! \right|<\infty, \left| \! \left| \! \left| w_j\right| \! \right| \! \right|<\infty, j\in \{1, 2\}$.
\end{lemma}
\begin{proof}
We first consider the $L^3$-norm of $\rho^{\sigma}H[w_1, w_2](t)$. Note that $\left| \! \left| \! \left| w_1\right| \! \right| \! \right|_{\sigma, t}<\infty$ implies $\rho^{\sigma}w_1(t)\in L^6(D)$ and 
\begin{align*}
\|\rho^{\sigma}w_1(t)\|_6\le t^{-\frac{1}{4}}(1+t)^{\sigma}[w_1]^{\frac{1}{2}}_{3, \sigma, t}[w_1]^{\frac{1}{2}}_{\infty, \sigma, t}.
\end{align*}
It follows from \eqref{FK-wei}, \eqref{dec-0-alpha1} and \eqref{dec-1-alpha1} that 
\begin{align*}
&\|\rho^{\sigma}H[w_1, w_2](t)\|_3\\
\le &C\int_{0}^{t}(t-\tau)^{-\frac{1}{4}}\left\{\|\rho^{\sigma}P(w_1\cdot\nabla w_2)(\tau)\|_{2}+(1+t-\tau)^{\sigma}\|P(w_1\cdot\nabla w_2)(\tau)\|_{2}\right\}\,d\tau\\
\le &C\int_{0}^{t}(t-\tau)^{-\frac{1}{4}}\left\{\|\rho^{\sigma}w_1(\tau)\|_6\|\nabla w_2(\tau)\|_3+(1+t-\tau)^{\sigma}\|w_1(\tau)\|_6\|\nabla w_2(\tau)\|_3\right\}\,d\tau\\
\le &C\int_{0}^{t}(t-\tau)^{-\frac{1}{4}}\tau^{-\frac{3}{4}}(1+\tau)^{\sigma}\,d\tau\;[w_1]^{\frac{1}{2}}_{3, \sigma, t}[w_1]^{\frac{1}{2}}_{\infty, \sigma,  t}[\nabla w_2]_{3, 0, t}\\
+ &C\int_{0}^{t}(t-\tau)^{-\frac{1}{4}}\tau^{-\frac{3}{4}}(1+t-\tau)^{\sigma}\,d\tau\;[w_1]^{\frac{1}{2}}_{3, 0, t}[w_1]^{\frac{1}{2}}_{\infty, 0,  t}[\nabla w_2]_{3, 0, t}\\
\le &C(1+t)^{\sigma}\left([w_1]^{\frac{1}{2}}_{3, 0, t}[w_1]^{\frac{1}{2}}_{\infty, 0, t}+[w_1]^{\frac{1}{2}}_{3, \sigma, t}[w_1]^{\frac{1}{2}}_{\infty, \sigma,  t}\right)[\nabla w_2]_{3, 0, t}
\end{align*}
and that
\begin{align*}
&\|\rho^{\sigma} H[w_1, w_2](t)\|_{\infty}+\|\rho^{\sigma}\nabla H[w_1, w_2](t)\|_{3}\\
\le &C\int_{0}^{t}(t-\tau)^{-\frac{3}{4}}\left\{\|\rho^{\sigma}P(w_1\cdot\nabla w_2)(\tau)\|_{2}+(1+t-\tau)^{\sigma}\|P(w_1\cdot\nabla w_2)(\tau)\|_{2}\right\}\,d\tau\\
\le &C\int_{0}^{t}(t-\tau)^{-\frac{3}{4}}\tau^{-\frac{3}{4}}(1+\tau)^{\sigma}\,d\tau [w_1]^{\frac{1}{2}}_{3, \sigma, t}[w_1]^{\frac{1}{2}}_{\infty, \sigma,  t}[\nabla w_2]_{3, 0, t}\\
+ &C\int_{0}^{t}(t-\tau)^{-\frac{3}{4}}\tau^{-\frac{3}{4}}(1+t-\tau)^{\sigma}\,d\tau\;[w_1]^{\frac{1}{2}}_{3, 0, t}[w_1]^{\frac{1}{2}}_{\infty, 0,  t}[\nabla w_2]_{3, 0, t}\\
\le &Ct^{-\frac{1}{2}}(1+t)^{\sigma}\left([w_1]^{\frac{1}{2}}_{3, 0, t}[w_1]^{\frac{1}{2}}_{\infty, 0, t}+[w_1]^{\frac{1}{2}}_{3, \sigma, t}[w_1]^{\frac{1}{2}}_{\infty, \sigma,  t}\right)[\nabla w_2]_{3, 0, t},
\end{align*}
with some constants $C>0$. Collecting these estimates, we have \eqref{hest}.\par
Let us consider the estimate of $G[w]$. 
We see from \eqref{dec-0-alpha1} that
\begin{align*}
\|\rho^{\sigma}G[w](t)\|_3
&\le C\int_{0}^{t}(t-\tau)^{-\frac{1}{2}}\left\{\|\rho^{\sigma}P(v\cdot\nabla w)(\tau)\|_{\frac{3}{2}}+(1+t-\tau)^{\sigma}\|P(v\cdot\nabla w)\|_{\frac{3}{2}}\right\}\,d\tau\\
&\le C\int_{0}^{t}(t-\tau)^{-\frac{1}{2}}\left\{\|v(\tau)\|_3\|\rho^{\sigma}\nabla w(\tau)\|_3+(1+t-\tau)^{\sigma}\|v(\tau)\|_3\|\nabla w(\tau)\|_3\right\}\,d\tau\\
&\le C (1+t)^{\sigma}\llbracket v\rrbracket \bigl([\nabla w]_{3, 0, t}+[\nabla w]_{3, \sigma, t}\bigr),
\end{align*}
which leads to $[G[w]]_{3, \sigma, t}\le C\llbracket v\rrbracket \left([\nabla w]_{3, 0, t}+[\nabla w]_{3, \sigma, t}\right)$.\par
To estimate $\|\rho^{\sigma}G[w](t)\|_{\infty}$ and $\|\nabla \rho^{\sigma}G[w](t)\|_{3}$, we will combine two estimates: $t\le 2$ and $t\ge 2$.
For $t\le 2$, \eqref{wei-est1} leads us to
\begin{align*}
&\|\rho^{\sigma}G[w](t)\|_{\infty}+\|\rho^{\sigma}\nabla G[w](t)\|_{3}\\
\le &C\int_{0}^{t}\left\{(t-\tau)^{-\frac{1}{2}}\|\rho^{\sigma}P(v\cdot\nabla w)(\tau)\|_{3}+(t-\tau)^{-\frac{1}{2}}\|\rho^{\sigma}P(w\cdot\nabla v)(\tau)\|_{3}\right\}\,d\tau\\
\le &C\int_{0}^{t}\left\{(t-\tau)^{-\frac{1}{2}}\|v(\tau)\|_{\infty}\|\rho^{\sigma}\nabla w(\tau)\|_3+(t-\tau)^{-\frac{1}{2}}\|\rho^{\sigma}w(\tau)\|_{\infty}\|\nabla v(\tau)\|_{3}\right\}\,d\tau\\
\le &C(1+t)^{\sigma}\llbracket v\rrbracket \bigl([w]_{\infty, \sigma, t}+[\nabla w]_{3, \sigma, t}\bigr),
\end{align*}
which implies $[G[w]]_{\infty, \sigma, t}+[\nabla G[w]]_{3, \sigma, t}\le C\llbracket v\rrbracket\bigl([w]_{\infty, \sigma, t}+[\nabla w]_{3, \sigma, t}\bigr)$.

For $t\ge 2$, we split the integral into three parts:
\begin{align}\label{g-split}
G[w](t)=\left(\int_0^{\frac{t}{2}}+\int_{\frac{t}{2}}^{t-1}+\int_{t-1}^{t}\right)=:G_1[w](t)+G_2[w](t)+G_3[w](t).
\end{align}
For $G_1[w]$, it follows from \eqref{dec-0-alpha1}--\eqref{dec-1-alpha1} that
\begin{align*}
&\|\rho^{\sigma}G_1[w](t)\|_{\infty}+\|\rho^{\sigma}\nabla G_1[w](t)\|_3\\
\le &C\int_{0}^{\frac{t}{2}}(t-\tau)^{-1-\frac{\delta_1}{2}}\left\{\|\rho^{\sigma}P(v\cdot\nabla w)(\tau)\|_{\frac{3}{2+\delta_1}}+(1+t-\tau)^{\sigma}\|P(v\cdot\nabla w)(\tau)\|_{\frac{3}{2+\delta_1}} \right\} \,d\tau\\
+ &C\int_{0}^{\frac{t}{2}}(t-\tau)^{-1-\frac{\delta_2}{2}}\left\{\|\rho^{\sigma}P(w\cdot\nabla v)(\tau)\|_{\frac{3}{2+\delta_2}}+(1+t-\tau)^{\sigma}\|P(w\cdot\nabla v)(\tau)\|_{\frac{3}{2+\delta_2}} \right\} \,d\tau\\
\le &C\int_{0}^{\frac{t}{2}}(t-\tau)^{-1-\frac{\delta_1}{2}}\left\{\|v(\tau)\|_{\frac{3}{1+\delta_1}}\|\rho^{\sigma}\nabla w(\tau)\|_3+(1+t-\tau)^{\sigma}\|v(\tau)\|_{\frac{3}{1+\delta_1}}\|\nabla w(\tau)\|_3 \right\}\,d\tau\\
+ &C\int_{0}^{\frac{t}{2}}(t-\tau)^{-1-\frac{\delta_2}{2}}\left\{\|\rho^{\sigma}w(\tau)\|_{\infty}\|\nabla v(\tau)\|_{\frac{3}{2+\delta_2}}+(1+t-\tau)^{\sigma}\|w(\tau)\|_{\infty}\|\nabla v(\tau)\|_{\frac{3}{2+\delta_2}} \right\}\,d\tau\\
\le &C\left(\frac{t}{2}\right)^{-1-\frac{\delta_1}{2}}\left\{\int_{0}^{\frac{t}{2}}\tau^{-\frac{1}{2}}(1+\tau)^{\sigma}\,d\tau\llbracket v\rrbracket [\nabla w]_{3, \sigma, t}+\int_{0}^{\frac{t}{2}}(1+t-\tau)^{\sigma}\tau^{-\frac{1}{2}}\,d\tau \llbracket v\rrbracket [\nabla w]_{3, 0, t}\right\}\\
+ &C\left(\frac{t}{2}\right)^{-1-\frac{\delta_2}{2}}\left\{\int_{0}^{\frac{t}{2}}\tau^{-\frac{1}{2}}(1+\tau)^{\sigma}\,d\tau\llbracket v\rrbracket [w]_{\infty, \sigma, t}+\int_{0}^{\frac{t}{2}}(1+t-\tau)^{\sigma}\tau^{-\frac{1}{2}}\,d\tau \llbracket v\rrbracket [w]_{\infty, 0, t}\right\}\\
\le &C(1+t)^{\sigma}\llbracket v\rrbracket \left\{t^{-\frac{1}{2}-\frac{\delta_1}{2}} \bigl([\nabla w]_{3, 0, t}+[\nabla w]_{3, \sigma, t}\bigr)+t^{-\frac{1}{2}-\frac{\delta_2}{2}}\bigl([w]_{\infty, 0, t}+[w]_{\infty, \sigma, t}\bigr)\right\}
\end{align*}
for $t\ge2$.
Remark that \eqref{def-delta} guarantees
\begin{align*}
\frac{3}{3-\alpha}<\frac{3}{2+\delta_1},\qquad \frac{3}{3-\alpha}<\frac{3}{2+\delta_2},
\end{align*}
and, therefore, the inequalities 
\begin{align*}
\|\rho^{\alpha}P(v\cdot\nabla w)\|_{\frac{3}{2+\delta_1}}\le C\|\rho^{\alpha}(v\cdot\nabla w)\|_{\frac{3}{2+\delta_1}},\qquad \|\rho^{\alpha}P(w\cdot\nabla v)\|_{\frac{3}{2+\delta_2}}\le C\|\rho^{\alpha}(w\cdot\nabla v)\|_{\frac{3}{2+\delta_2}}
\end{align*}
hold true. 
As for $G_2[w]$, we employ \eqref{dec-0-alpha1}--\eqref{dec-1-alpha1} to find
\begin{align*}
&\|\rho^{\sigma}G_2[w](t)\|_{\infty}+\|\rho^{\sigma}\nabla G_2[w](t)\|_{3}\\
\le &C\int_{\frac{t}{2}}^{t-1}(t-\tau)^{-1-\frac{\delta_1}{2}}\left\{\|v(\tau)\|_{\frac{3}{1+\delta_1}}\|\rho^{\sigma}\nabla w(\tau)\|_3+(1+t-\tau)^{\sigma}\|v(\tau)\|_{\frac{3}{1+\delta_1}}\|\nabla w(\tau)\|_3 \right\}\,d\tau\\
+ &C\int_{\frac{t}{2}}^{t-1}(t-\tau)^{-1-\frac{\delta_2}{2}}\left\{\|\rho^{\sigma}w(\tau)\|_{\infty}\|\nabla v(\tau)\|_{\frac{3}{2+\delta_2}}+(1+t-\tau)^{\sigma}\|w(\tau)\|_{\infty}\|\nabla v(\tau)\|_{\frac{3}{2+\delta_2}} \right\}\,d\tau\\
\le &C\int_{\frac{t}{2}}^{t-1}(t-\tau)^{-1-\frac{\delta_1}{2}}\tau^{-\frac{1}{2}}(1+\tau)^{\sigma}\,d\tau\llbracket v\rrbracket[\nabla w]_{3, \sigma, t}+C\int_{\frac{t}{2}}^{t-1}(t-\tau)^{-1-\frac{\delta_1}{2}}(1+t-\tau)^{\sigma}\tau^{-\frac{1}{2}}\,d\tau\llbracket v\rrbracket[\nabla w]_{3, 0, t}\\
+ &C\int_{\frac{t}{2}}^{t-1}(t-\tau)^{-1-\frac{\delta_2}{2}}\tau^{-\frac{1}{2}}(1+\tau)^{\sigma}\,d\tau\llbracket v\rrbracket [w]_{\infty, \sigma, t}+\int_{\frac{t}{2}}^{t-1}(t-\tau)^{-1-\frac{\delta_2}{2}}(1+t-\tau)^{\sigma}\,d\tau\llbracket v\rrbracket [w]_{\infty, 0, t}\\
\le &C(1+t)^{\sigma}\left(\frac{t}{2}\right)^{-\frac{1}{2}}\llbracket v\rrbracket \bigl([\nabla w_1]_{3, 0, t}+[\nabla w_1]_{3, \sigma, t}+[w]_{\infty, 0, t}+[w]_{\infty, \sigma, t}\bigr).
\end{align*}
For $G_3[w]$, estimate \eqref{est-uni-wei} implies that
\begin{align*}
&\|\rho^{\sigma}G_3[w](t)\|_{\infty}+\|\rho^{\sigma}\nabla G_3[w](t)\|_3\\
\le &C\int_{t-1}^{t}(t-\tau)^{-\frac{1}{2}}\|\rho^{\sigma}P(v\cdot\nabla w)(\tau)\|_{3}\,d\tau
+C\int_{t-1}^{t}(t-\tau)^{-\frac{1}{2}}\|\rho^{\sigma}P(w\cdot\nabla v)(\tau)\|_{3}\,d\tau\\
\le &C\int_{t-1}^{t}(t-\tau)^{-\frac{1}{2}}\|v(\tau)\|_{\infty}\|\rho^{\sigma}\nabla w(\tau)\|_{3}\,d\tau
+C\int_{t-1}^{t}(t-\tau)^{-\frac{1}{2}}\|\rho^{\sigma}w(\tau)\|_{\infty}\|\nabla v(\tau)\|_{3}\,d\tau\\
\le &C\int_{t-1}^{t}(t-\tau)^{-\frac{1}{2}}\tau^{-\frac{1}{2}}(1+\tau)^{\sigma}\,d\tau\llbracket v\rrbracket[\nabla w]_{3, \sigma, t}
+C\int_{t-1}^{t}(t-\tau)^{-\frac{1}{2}}\tau^{-\frac{1}{2}}(1+\tau)^{\sigma}\,d\tau\llbracket v\rrbracket[w]_{\infty, \sigma, t}\\
\le &C(1+t)^{\sigma}t^{-\frac{1}{2}}\llbracket v\rrbracket\bigl([\nabla w]_{3, \sigma, t}+[w]_{\infty, \sigma, t}\bigr).
\end{align*}
Collecting these estimates leads to \eqref{gest}.
The proof is complete. 
\end{proof}

For later use, we set constants
\begin{align*}
C_1^{\prime}:=2\mathcal{C}_{0},\qquad 
C_1^{\prime\prime}:=\mathcal{C}_{\alpha},\qquad 
C_2^{\prime}:=2\widetilde{\mathcal{C}}_{0},\qquad
C_2^{\prime\prime}:=\widetilde{\mathcal{C}}_{\alpha},\qquad
C_1:=C_1^{\prime}+C_1^{\prime\prime},\qquad
C_2:=C_2^{\prime}+C_2^{\prime\prime},
\end{align*}
where $\mathcal C_{\sigma}$ and $\widetilde{\mathcal C}_{\sigma}$ with $\sigma\in\{0, \alpha\}$ are the constants in \eqref{hest}--\eqref{gest}.
Then we see from Lemma \ref{stab-lem-duha} that
\begin{align}
&\left| \! \left| \! \left| H[w_1, w_2]\right| \! \right| \! \right|_{0, t}\le C_1^{\prime}[w_1]^{\frac{1}{2}}_{3, 0, t}[w_1]^{\frac{1}{2}}_{\infty, 0, t}[\nabla w_2]_{3, 0, t},\label{h-0-const}\\
&\left| \! \left| \! \left| G[w]\right| \! \right| \! \right|_{0, t}\le C_2^{\prime}\llbracket v\rrbracket \bigl([\nabla w]_{3, 0, t}+[w]_{\infty, 0, t}\bigr)\label{g-0-const},\\
&\left| \! \left| \! \left| H[w_1, w_2]\right| \! \right| \! \right|_{\alpha, t}\le C_1^{\prime\prime}\left([w_1]^{\frac{1}{2}}_{3, 0, t}[w_1]^{\frac{1}{2}}_{\infty, 0, t}+[w_1]^{\frac{1}{2}}_{3, \alpha, t}[w_1]^{\frac{1}{2}}_{\infty, \alpha, t}\right)[\nabla w_2]_{3, 0, t},\label{h-wei-const}\\
&\left| \! \left| \! \left| G[w]\right| \! \right| \! \right|_{\alpha, t}\le C_2^{\prime\prime}\llbracket v\rrbracket \bigl([\nabla w]_{3, 0, t}+[\nabla w]_{3, \alpha, t}+[w]_{\infty, 0, t}+[w]_{\infty, \alpha, t}\bigr).\label{g-wei-const}\\
&\left| \! \left| \! \left| H[w_1, w_2]\right| \! \right| \! \right|_{0, t}+\left| \! \left| \! \left| H[w_1, w_2]\right| \! \right| \! \right|_{\alpha, t}\le C_1\left([w_1]^{\frac{1}{2}}_{3, 0, t}[w_1]^{\frac{1}{2}}_{\infty, 0, t}+[w_1]^{\frac{1}{2}}_{3, \alpha, t}[w_1]^{\frac{1}{2}}_{\infty, \alpha,  t}\right)[\nabla w_2]_{3, 0, t}\label{h-x-const},\\
&\left| \! \left| \! \left| G[w]\right| \! \right| \! \right|_{0, t}+\left| \! \left| \! \left| G[w]\right| \! \right| \! \right|_{\alpha, t}\le C_2\llbracket v\rrbracket \bigl([\nabla w]_{3, 0, t}+[\nabla w]_{3, \alpha, t}+[w]_{\infty, 0, t}+[w]_{\infty, \alpha, t}\bigr)\label{g-x-const}.
\end{align}

\begin{remark}\label{rmk-const-j}
Let $\mathcal T_0>0$ and $t\ge \mathcal T_0$.
We consider the integrals
\begin{gather}
\widetilde{H}[w_1, w_2](t):=\int_{\mathcal T_0}^{t}T(t, \tau)P(w_1\cdot\nabla w_2)(\tau)\,d\tau,\label{int-h-t0}\\
\widetilde{G}[w](t):=\int_{\mathcal T_0}^{t}T(t, \tau)P(v\cdot\nabla w+w\cdot\nabla v)(\tau)\,d\tau,\label{int-g-t0}
\end{gather}
and set
\begin{gather}
[w]_{q, \sigma, J}:=\sup_{\tau\in J}(\tau-\mathcal T_0)^{\frac{1}{2}-\frac{3}{2q}}(1+\tau)^{-\sigma}\|\rho^{\sigma}w(\tau)\|_q,\label{sob-infty-t0}\\
[\nabla w]_{3, \sigma, J}:=\sup_{\tau\in J}(\tau-\mathcal T_0)^{\frac{1}{2}}(1+\tau)^{-\sigma}\|\rho^{\sigma}\nabla w(\tau)\|_3,\label{sob-nab-t0}
\end{gather}
for $\sigma\in \{0, \alpha\}$, where $J:=[\mathcal T_0, t]$.
The same computations as in Lemma \ref{stab-lem-duha}, in which the formula \eqref{g-split} is replaced by
\begin{align*}
\int_{\mathcal T_0}^{t}=\int_{\mathcal T_0}^{\frac{t+\mathcal T_0}{2}}+\int_{\frac{t+\mathcal T_0}{2}}^{t-1}+\int_{t-1}^{t}
\end{align*}
for $t-\mathcal T_0\ge2$, lead us to
\begin{align}
&\left[\widetilde{H}[w_1, w_2]\right]_{\infty, 0, J}+\left[\nabla \widetilde{H}[w_1, w_2]\right]_{3, 0, J}\le C_1^{\prime}[w_1]^{\frac{1}{2}}_{3, 0, J}[w_1]^{\frac{1}{2}}_{\infty, 0, J}[\nabla w_2]_{3, 0, J},\label{h-0-j-const}\\
&\left[\widetilde{H}[w_1, w_2]\right]_{\infty, \alpha, J}+\left[\nabla \widetilde{H}[w_1, w_2]\right]_{3, \alpha, J}\le C_1^{\prime\prime}\left([w_1]^{\frac{1}{2}}_{3, 0, J}[w_1]^{\frac{1}{2}}_{\infty, 0, J}+[w_1]^{\frac{1}{2}}_{3, \alpha, J}[w_1]^{\frac{1}{2}}_{\infty, \alpha, J}\right)[\nabla w_2]_{3, 0, J},\label{h-wei-j-const}\\
&\left[\widetilde{G}[w]\right]_{\infty, 0, J}+\left[\nabla\widetilde{G}[w]\right]_{3, 0, J}\le C_2^{\prime}\llbracket v\rrbracket \bigl([\nabla w]_{3, 0, J}+[w]_{\infty, 0, J}\bigr),\label{g-0-j-const}\\
&\left[\widetilde{G}[w]\right]_{\infty, \alpha, J}+\left[\nabla\widetilde{G}[w]\right]_{3, \alpha, J}\le C_2^{\prime\prime}\llbracket v\rrbracket \bigl([\nabla w]_{3, 0, J}+[\nabla w]_{3, \alpha, J}+[w]_{\infty, 0, J}+[w]_{\infty, \alpha, J}\bigr),\label{g-wei-j-const}
\end{align}
where the constants $C_1^{\prime}, C_1^{\prime\prime}, C_2^{\prime}, C_2^{\prime\prime}$ are as in \eqref{h-0-const}--\eqref{g-wei-const}.
Estimates \eqref{h-0-j-const}--\eqref{g-wei-j-const} will be utilized to show the sharp behavior \eqref{sol-0-sob}--\eqref{sol-wei-sob}.
\end{remark}

Based on \eqref{h-x-const}--\eqref{g-x-const}, it is possible to show the existence of a solution to \eqref{inteq}.
To this end, we introduce the Banach space
\begin{align}\label{stab-def-ball}
X:=\{w\in L^{\infty}(0, \infty; L^3_{\sigma}(D));\; \| w\|_X<\infty,\; \lim_{t\to0}\left(\|\rho^{\alpha}w(t)-\rho^{\alpha}w_0\|_3+\|w\|_{X(t)}\right)=0\},
\end{align}
where
\begin{align*}
\|w\|_{X}:=\sum_{\sigma=0, \alpha}\sup_{t> 0}\left| \! \left| \! \left| w\right| \! \right| \! \right|_{\sigma, t},\qquad \| w\|_{X(t)}:=\sum_{\sigma=0, \alpha}\left([w]_{\infty, \sigma, t}+[\nabla w]_{3, \sigma, t}\right),
\end{align*}
and set $X_R:=\{w\in X;\; \|w\|_X\le R\}$ with radius $R>0$ to be determined later.

\begin{proposition}\label{stab-prop-exist}
In addition to the assumptions of Proposition \ref{prop-galtak}, let $\alpha\in [0, 1)$ and $w_0\in L^3_{\alpha, \sigma}(D)$. There are constants $\varepsilon_0, \varepsilon_1>0$ such that if
\begin{align*}
\|\rho^{\alpha}w_0\|_3\le \varepsilon_0,\qquad \lambda^{-\frac{3}{4}}\|(\eta, \omega)\|_{W^{2, 2}(0, \mathcal T)}\le \varepsilon_1,
\end{align*}
as well as \eqref{per-as-small}.
Then \eqref{inteq} admits a unique solution $w\in X$ with
\begin{align}\label{sol-rad}
\|w\|_{X}\le 4C_0\|\rho^{\alpha}w_0\|_3
\end{align}
and 
\begin{align*}
\lim_{t\to0}\|\rho^{\alpha}w(t)-\rho^{\alpha}w_0\|_3=0,
\end{align*}
where $C_0>0$ is given by \eqref{def-c0} below.
\end{proposition}

\begin{proof}
Define the map $\Phi$ by
\begin{align*}
\Phi[w](t)=T(t, 0)w_0-H[w, w](t)-G[w](t).
\end{align*}  
We find a unique solution to \eqref{inteq} through the Banach fixed point theorem. 
By \eqref{dec-0-alpha1} and \eqref{dec-1-alpha1}, there is a constant $C_0>0$ such that
\begin{align}\label{def-c0}
\| T(\cdot, 0)w_0\|_{X}\le C_0\|\rho^{\alpha}w_0\|_3.
\end{align}
We see from \eqref{wei-small1}--\eqref{IC-1} that
\begin{align*}
\|\rho^{\alpha}T(t, 0)w_0-\rho^{\alpha}w_0\|_3\to 0,\qquad \|T(\cdot, 0)w_0\|_{X(t)}\to 0
\end{align*}
as $t\to 0$.
As mentioned in Remark \ref{rmk-frac34}, there is a constant $c_*=c_*(\delta_1, \delta_2)>0$ such that
\begin{align}\label{cond-small-per}
\llbracket v \rrbracket \le c_*\lambda^{-\frac{3}{4}}\|(\eta, \omega)\|_{W^{2, 2}(0, \mathcal T)}.
\end{align}
From this together with \eqref{h-x-const}--\eqref{g-x-const}, we see that
\begin{gather*}
\|H[w, w]\|_{X}\le 2C_1R^2,\qquad \|G[w]\|_X\le 4C_2c_*\lambda^{-\frac{3}{4}}\|(\eta, \omega)\|_{W^{2, 2}(0, \mathcal T)} R,\\
\|H[w, w]\|_{X(t)}+\|G[w]\|_{X(t)}\to 0 \quad \mbox{as}\quad t\to0
\end{gather*}
for all $w\in X_R$ with constants $C_1, C_2>0$ given by \eqref{h-x-const}--\eqref{g-x-const}.
Moreover, it follows from Lemma \ref{rmk-frac34} again that
\begin{align*}
&\|H[w_1, w_1]-H[w_2, w_2]\|_X\le 4C_1R\|w_1-w_2\|_X,\\
&\|G[w_1]-G[w_2]\|_X\le 4C_2c_*\lambda^{-\frac{3}{4}}\|(\eta, \omega)\|_{W^{2, 2}(0, \mathcal T)} \|w_1-w_2\|_X
\end{align*}
for every $w_1, w_2\in X_R$. 
We thus obtain
\begin{align*}
&\|\Phi[w]\|_X
\le 2C_1R^2+4C_2c_*\lambda^{-\frac{3}{4}}\|(\eta, \omega)\|_{W^{2, 2}(0, \mathcal T)}R+C_0\|\rho^{\alpha}w_0\|_3\\
&\|\Phi[w_1]-\Phi[w_2]\|_X\le (4C_1R+4C_2c_*\lambda^{-\frac{3}{4}}\|(\eta, \omega)\|_{W^{2, 2}(0, \mathcal T)})\|w_1-w_2\|_X
\end{align*}
for $w, w_1, w_2\in X_R$. If we assume
\begin{align}\label{cond-small-ex}
\|\rho^{\alpha}w_0\|_{3}\le \frac{1}{32C_0C_1}=:\varepsilon_0, \qquad \lambda^{-\frac{3}{4}}\|(\eta, \omega)\|_{W^{2, 2}(0, \mathcal T)}\le \frac{1}{8C_2c_*}=:\varepsilon_1,
\end{align}
and set
\begin{align*}
R:= \frac{1}{4C_1}\left(\frac{1}{2}-\sqrt{\frac{1}{4}-8C_0C_1\|\rho^{\alpha}w_0}\|_3\right),
\end{align*}
then we have 
\begin{align*}
\|\Phi[w]\|_{X}\le R< 4C_0\|\rho^{\alpha}w_0\|_3
\end{align*}
as well as
\begin{align*}
4C_1R+4C_2c_*\lambda^{-\frac{3}{4}}\|(\eta, \omega)\|_{W^{2, 2}(0, \mathcal T)}<1.
\end{align*}
Therefore $\Phi$ is contractive from $X_{R}$ into itself with $R>0$ determined above. 
The proof is complete.
\end{proof}

To obtain the sharp asymptotic behavior \eqref{sol-wei-sob}, the following lemma is the first step.\\
\begin{lemma}\label{stab-lem-sob-app}
In addition to the assumptions of Proposition \ref{stab-prop-exist}, suppose $w_0\in C^{\infty}_{0, \sigma}(D)$.
There exists a constant $\varepsilon_2^{\prime}\in (0, \varepsilon_0]$  such that if 
\begin{align}\label{sob-approx}
\|\rho^{\alpha}w_0\|_3\le \varepsilon_2^{\prime},
\end{align}
then the solution $w\in X$ to \eqref{inteq} obtained in Proposition \ref{stab-prop-exist} enjoys 
\begin{align}\label{aux-w-sob-app}
\|w(t)\|_3=o(1),\qquad \|\rho^{\alpha}w(t)\|_3=o\left(t^{\alpha}\right)
\end{align}
as $t\to\infty$.
\end{lemma}

\begin{proof}
We set $\delta:=\min\{\delta_1, \delta_2\}$ where $\delta_1$ and $\delta_2$ are as in \eqref{def-delta}.
By Theorem \ref{thm1}, we have
\begin{align}
t^{\frac{\delta}{2}}\|T(t, 0)w_0\|_3\le c_0\|w_0\|_{\frac{3}{1+\delta}},\qquad 
t^{\frac{\delta}{2}}(1+t)^{-\alpha}\|\rho^{\alpha}T(t, 0)w_0\|_{\frac{3}{1+\delta}}
\le c^{\prime}_{0}\|\rho^{\alpha}w_0\|_{\frac{3}{1+\delta}}\label{sob-top-delta}
\end{align}
with some constants $c_{0}, c^{\prime}_{0}>0$.
Let $\sigma\in \{0, \alpha\}$. We see from \eqref{dec-0-alpha1} that
\begin{align}
&\|\rho^{\sigma}H[w, w](t)\|_3\notag\\
\le &C\int_{0}^{t}(t-\tau)^{-\frac{1}{2}}\left\{\|\rho^{\sigma}P(w\cdot\nabla w)(\tau)\|_{\frac{3}{2}}+(1+t-\tau)^{\sigma}\|P(w\cdot\nabla w)(\tau)\|_{\frac{3}{2}}\right\}\,d\tau\notag\\
\le &C\int_{0}^{t}(t-\tau)^{-\frac{1}{2}}\tau^{-\frac{1}{2}-\frac{\delta}{2}}(1+\tau)^{\sigma}\,d\tau\;[\nabla w]_{3, 0, t}\left(\sup_{0<\tau<t}\tau^{\frac{\delta}{2}}(1+\tau)^{-\sigma}\|\rho^{\sigma}w(\tau)\|_3\right)\notag\\
+ &C\int_{0}^{t}(t-\tau)^{-\frac{1}{2}}(1+t-\tau)^{\sigma}\tau^{-\frac{1}{2}-\frac{\delta}{2}}\,d\tau\;[\nabla w]_{3, 0, t}\left(\sup_{0<\tau<t}\tau^{\frac{\delta}{2}}\|w(\tau)\|_3\right)\notag\\
\le &Ct^{-\frac{\delta}{2}}(1+t)^{\sigma}\|w\|_X\left(\sup_{0<\tau<t}\tau^{\frac{\delta}{2}}(1+\tau)^{-\sigma}\|\rho^{\sigma}w(\tau)\|_3+\sup_{0<\tau<t}\tau^{\frac{\delta}{2}}\|w(\tau)\|_3\right).\label{sob-nonlin-delta}
\end{align}
As for $\|\rho^{\sigma}G[w](t)\|_3$, it follows from \eqref{wei-est1} with $j=0$ that
\begin{align}
\|\rho^{\sigma}G[w](t)\|_3
\le &C\int_{0}^{t}(t-\tau)^{-\frac{1}{2}}\|\rho^{\sigma}P(v\cdot\nabla w+w\cdot\nabla v)(\tau)\|_{\frac{3}{2}}\,d\tau\notag\\
\le &C\int_{0}^{t}(t-\tau)^{-\frac{1}{2}}\left(\|v(\tau)\|_{3}\|\rho^{\sigma}\nabla w(\tau)\|_3+\|\rho^{\sigma}w(\tau)\|_{\infty}\|\nabla v(\tau)\|_3\right)\,d\tau\notag\\
\le &C\int_{0}^{t}(t-\tau)^{-\frac{1}{2}}\tau^{-\frac{1}{2}}(1+\tau)^{\sigma}\,d\tau\llbracket v\rrbracket \left([\nabla w]_{3, \sigma, t}+[w]_{\infty, \sigma, t}\right)\notag\\
\le &C(1+t)^{\sigma}\llbracket v\rrbracket \|w\|_{X}\label{sob-g-short}
\end{align}
for $t\le 2$.
On the other hand, estimate \eqref{dec-0-alpha1} leads us to
\begin{align}
&\|\rho^{\sigma}G[w](t)\|_3\notag\\
\le &C\int_{0}^{t}(t-\tau)^{-\frac{1}{2}-\frac{\delta_1}{2}}\left\{\|\rho^{\sigma}P(v\cdot\nabla w)(\tau)\|_{\frac{3}{2+\delta_1}}+(1+t-\tau)^{\sigma}\|P(v\cdot\nabla w)(\tau)\|_{\frac{3}{2+\delta_1}}\right\}\,d\tau \notag\\
+ &C\int_{0}^{t}(t-\tau)^{-\frac{1}{2}-\frac{\delta_2}{2}}\left\{\|\rho^{\sigma}P(w\cdot\nabla v)(\tau)\|_{\frac{3}{2+\delta_2}}+(1+t-\tau)^{\sigma}\|P(w\cdot\nabla v)(\tau)\|_{\frac{3}{2+\delta_2}}\right\}\,d\tau\notag\\
\le &C\int_{0}^{t}(t-\tau)^{-\frac{1}{2}-\frac{\delta_1}{2}}\tau^{-\frac{1}{2}}(1+\tau)^{\sigma}\,d\tau\,\llbracket v\rrbracket [\nabla w]_{3, \sigma, t}+C\int_{0}^{t}(t-\tau)^{-\frac{1}{2}-\frac{\delta_1}{2}}(1+t-\tau)^{\sigma}\tau^{-\frac{1}{2}}\,d\tau\,\llbracket v\rrbracket [\nabla w]_{3, 0, t}\notag\\
+ &C\int_{0}^{t}(t-\tau)^{-\frac{1}{2}-\frac{\delta_2}{2}}\tau^{-\frac{1}{2}}(1+\tau)^{\sigma}\,d\tau\,\llbracket v\rrbracket [w]_{\infty, \sigma, t}+C\int_{0}^{t}(t-\tau)^{-\frac{1}{2}-\frac{\delta_2}{2}}(1+t-\tau)^{\sigma}\tau^{-\frac{1}{2}}\,d\tau\,\llbracket v\rrbracket [w]_{\infty, 0, t}\notag\\
\le &Ct^{-\frac{\delta}{2}}(1+t)^{\sigma}\llbracket v\rrbracket \|w\|_{X}\label{sob-g-large}
\end{align}
for $t\ge2$. By \eqref{sob-g-short}--\eqref{sob-g-large}, we find
\begin{align}
\|\rho^{\sigma}G[w](t)\|_{3}\le C(1+t)^{\sigma-\frac{\delta}{2}}\llbracket v\rrbracket \|w\|_X,\label{sob-g-delta}
\end{align}
for all $t>0$.
Collecting estimates \eqref{sob-top-delta}--\eqref{sob-g-delta} leads to
\begin{align}
&\sup_{0<\tau<t}\tau^{\frac{\delta}{2}}\|w(\tau)\|_3\le c_0\|w_0\|_{\frac{3}{1+\delta}}+c_1\|w\|_X\left(\sup_{0<\tau<t}\tau^{\frac{\delta}{2}}\|w(\tau)\|_3\right)+c_2\llbracket v\rrbracket \|w\|_X,\label{sob-w-0}
\end{align}
and
\begin{align}
&\sup_{0<\tau<t}\tau^{\frac{\delta}{2}}(1+t)^{-\alpha}\|\rho^{\alpha}w(\tau)\|_3\notag\\
\le &c^{\prime}_0\|\rho^{\alpha}w_0\|_{\frac{3}{1+\delta}}+c^{\prime}_1\|w\|_X\left(\sup_{0<\tau<t}\tau^{\frac{\delta}{2}}(1+t)^{-\alpha}\|\rho^{\alpha}w(\tau)\|_3+\sup_{0<\tau<t}\tau^{\frac{\delta}{2}}\|w(\tau)\|_3\right)+c^{\prime}_2\llbracket v\rrbracket \|w\|_X,\label{sob-w-wei}
\end{align}
where the constants $c_1, c_1^{\prime}, c_2$ and $c_2^{\prime}$ are dependent on $\alpha$ but independent of $t$.
Assume that
\begin{align*}
\|\rho^{\alpha}w_0\|_3\le\min\left\{\varepsilon_0, \frac{1}{8C_0c_1}\right\}=:\widetilde{\varepsilon_2^{\prime}},
\end{align*}
then it follows from \eqref{sol-rad} that $c_1\|w\|_X\le 1/2$, which combined with \eqref{sob-w-0} implies
\begin{align}\label{sob-3-app}
\sup_{0<\tau<t}\tau^{\frac{\delta}{2}}\|w(\tau)\|_3
&\le 2c_0\|w_0\|_{\frac{3}{1+\delta}}+8C_0c_2\llbracket v\rrbracket \|\rho^{\alpha}w_0\|_3,
\end{align}
thus, we have $\|w(t)\|_{3}=o(1)$.
Moreover, suppose
\begin{align}\label{sob-3-small}
\|\rho^{\alpha}w_0\|_3\le \min\left\{\widetilde{\varepsilon_2^{\prime}}, \frac{1}{8C_0c_1^{\prime}}\right\}=:\varepsilon_2^{\prime},
\end{align}
then we find $c_1^{\prime}\|w\|_X\le 1/2$ and \eqref{sob-3-app}.
Thus, we see from \eqref{sob-w-wei} that 
\begin{align*}
\sup_{0<\tau<t}\tau^{\frac{\delta}{2}}(1+\tau)^{-\alpha}\|\rho^{\alpha}w(\tau)\|_3
&\le 2c_0^{\prime}\|\rho^{\alpha}w_0\|_{\frac{3}{1+\delta}}+8C_0c_1^{\prime}\|\rho^{\alpha}w_0\|_3\left(2c_0\|w_0\|_{\frac{3}{1+\delta}}+8C_0c_2\llbracket v\rrbracket \|\rho^{\alpha}w_0\|_3\right)\\
&+8C_0c_2^{\prime}\|\rho^{\alpha}w_0\|_3\llbracket v\rrbracket
\end{align*}
for all $t>0$, yielding \eqref{aux-w-sob-app}.
The proof is complete.
\end{proof}

Let us close the paper with completion of the proof of Theorem \ref{thm2}.

\medskip
\noindent
{\it Proof of Theorem \ref{thm2}.}
Let $w$ be the solution to \eqref{inteq} obtained in Proposition \ref{stab-prop-exist}.
The regularity properties \eqref{sol-conti-l3}--\eqref{sol-conti-1st} were shown in \cite[Theorem 2.9]{T26}, and \eqref{sol-wei} is obvious due to the definiton of the function space $X$, see \eqref{stab-def-ball}.
Uniqueness within the class \eqref{uni-tak} was already clarified in \cite[Proposition 5.12]{T26}.
\par
Here, we will show \eqref{sol-0-sob}--\eqref{sol-wei-sob}.
We first consider \eqref{sol-0-sob} with $q=3$.
Let $w_0\in L^3_{\alpha, \sigma}(D)$ be a general initial disturbance satisfying
\begin{align}\label{def-ep2pp}
\|\rho^{\alpha}w_0\|_3\le \varepsilon_2^{\prime\prime},\qquad \varepsilon_2^{\prime\prime}:=\frac{\varepsilon_2^{\prime}}{2},
\end{align}
where $\varepsilon_2^{\prime}$ is given by \eqref{sob-approx}.
For any $\iota\in (0, \varepsilon_2^{\prime\prime})$, we choose $\overline{w_0}\in C^{\infty}_{0, \sigma}(D)$ such that $\|\rho^{\alpha}w_0-\rho^{\alpha}\overline{w_0}\|_3\le \iota$.
It follows from the triangle inequality that $\|\rho^{\alpha}\overline{w_0}\|_{3}\le \varepsilon_2^{\prime}$, thus, by Lemma \ref{stab-lem-sob-app}, we can find a solution $\overline{w}\in X$ to \eqref{inteq} with the initial disturbance $\overline{w_0}$ fulfilling \eqref{aux-w-sob-app}. 
We consider the difference
\begin{align*}
w(t)-\overline{w}(t)=T(t, 0)(w_0-\overline{w_0})-H[w-\overline{w}, w]-H[\overline{w}, w-\overline{w}]-G[w-\overline{w}]
\end{align*}
where $w\in X$ is the solution to \eqref{inteq} when the initial disturbance is $w_0\in L^3_{\alpha, \sigma}(D)$.
We immediately know from Theorem \ref{thm1} that
\begin{align}\label{sob-ini-app}
\|T(\cdot, 0)(w_0-\overline{w_0})\|_X
\le C_0\|\rho^{\alpha}w_0-\rho^{\alpha}\overline{w_0}\|_3.
\end{align}
From Lemma \ref{stab-lem-duha}, we have
\begin{align}
\| H[w-\overline{w}, w]-H[\overline{w}, w-\overline{w}]\|_X
&\le 2C_1\left(\|w\|_X+\|\overline{w}\|_X\right)\|w-\overline{w}\|_X\notag\\
&\le 8C_0C_1\left(\|\rho^{\alpha}w_0\|_3+\|\rho^{\alpha}\overline{w_0}\|_3\right)\|w-\overline{w}\|_X,\label{sob-nonlin-app}
\end{align}
and
\begin{align}
\|G[w-\overline{w}]\|_X
&\le 4C_2\llbracket v\rrbracket \|w-\overline{w}\|_X 
\le 4C_2c_*\lambda^{-\frac{3}{4}}\|(\eta, \omega)\|_{W^{2, 2}(0, \mathcal T)}\|w-\overline{w}\|_X.\label{sob-per-app}
\end{align}
It follows from \eqref{cond-small-ex}, \eqref{sob-3-small} and \eqref{def-ep2pp} that
\begin{align*}
\|\rho^{\alpha}w_0\|_3+\|\rho^{\alpha}\overline{w_0}\|_3\le \frac{3}{64C_0C_1},\qquad 4C_2c_*\lambda^{-\frac{3}{4}}\|(\eta, \omega)\|_{W^{2, 2}(0, \mathcal T)}\le \frac{1}{2},
\end{align*}
then \eqref{sob-ini-app}--\eqref{sob-per-app} lead us to
\begin{align*}
\|w-\overline{w}\|_X\le 8C_0\|\rho^{\alpha}w_0-\rho^{\alpha}\overline{w_0}\|_3.
\end{align*}
Hence, we get
\begin{align*}
(1+t)^{-\sigma}\|\rho^{\sigma}w(t)-\rho^{\sigma}\overline{w}(t)\|_3\le 8C_0\|\rho^{\alpha}w_0-\rho^{\alpha}\overline{w_0}\|_3
\end{align*}
and
\begin{align*}
(1+t)^{-\sigma}\|\rho^{\sigma}w(t)\|_3
&\le (1+t)^{-\sigma}\|\rho^{\sigma}w(t)-\rho^{\sigma}\overline{w}(t)\|_3+(1+t)^{-\sigma}\|\rho^{\sigma}\overline{w}(t)\|_3\\
&\le 8C_0\|\rho^{\alpha}w_0-\rho^{\alpha}\overline{w_0}\|_3+(1+t)^{-\sigma}\|\rho^{\sigma}\overline{w}(t)\|_3
\end{align*}
for $\sigma\in\{0, \alpha\}$.
By \eqref{aux-w-sob-app}, we obtain
\begin{align*}
\varlimsup_{t\to\infty}\|w(t)\|_3\le 8C_0\iota,\qquad \varlimsup_{t\to\infty}(1+t)^{-\alpha}\|\rho^{\alpha}w(t)\|_3\le 8C_0\iota.
\end{align*}
Since $\iota\in (0, \varepsilon_2^{\prime\prime})$ is arbitrary, we are led to \eqref{sol-0-sob}--\eqref{sol-wei-sob} with $q=3$.
By the interpolation estimate, we immediately have
\begin{align*}
&t^{\frac{1}{2}-\frac{3}{2q}}\|w(t)\|_q\le \left(\|w(t)\|_{3}\right)^{\frac{3}{q}}\left(t^{\frac{1}{2}}\|w(t)\|_{\infty}\right)^{1-\frac{3}{q}}\le (\|w(t)\|_{3})^{\frac{3}{q}}\|w\|_X^{1-\frac{3}{q}}\to 0,\\
&t^{\frac{1}{2}-\frac{3}{2q}-\alpha}\|\rho^{\alpha}w(t)\|_{q}\le \left(t^{-\alpha}\|\rho^{\alpha}w(t)\|_3\right)^{\frac{3}{q}}\left(t^{\frac{1}{2}-\alpha}\|\rho^{\alpha}w(t)\|_{\infty}\right)^{1-\frac{3}{q}}\le (t^{-\alpha}\|\rho^{\alpha}w(t)\|_3)^{\frac{3}{q}}\|w\|_X^{1-\frac{3}{q}}\to0
\end{align*}
as $t\to\infty$ for $q\in (3, \infty)$.
\par
We next prove 
\begin{align}\label{sob-step2}
\|w(t)\|_{\infty}+\|\nabla w(t)\|_3=o\left(t^{-\frac{1}{2}}\right),\qquad \|\rho^{\alpha}w(t)\|_{\infty}+\|\rho^{\alpha}\nabla w(t)\|_{3}=o\left(t^{-\frac{1}{2}+\alpha}\right)
\end{align}
by following the same argument as in \cite{T24}.
Let $\mathcal T_0>0$ and $t\ge \mathcal T_0$.
The integral equation \eqref{inteq} can be converted into 
\begin{align*}
w(t)=T(t, \mathcal T_0)w(\mathcal T_0)-\widetilde{H}[w, w](t)-\widetilde{G}[w](t),
\end{align*}
where $\widetilde{H}$ and $\widetilde{G}$ are given by \eqref{int-h-t0}--\eqref{int-g-t0}.
It follows from \eqref{lqlr-hi-0}--\eqref{lqlr-hi-1} that
\begin{align}\label{sob-0-top}
\sup_{\tau\in J}(\tau-\mathcal T_0)^{\frac{1}{2}}\|T(\tau, \mathcal T_0)w(\mathcal T_0)\|_{\infty}+\sup_{\tau\in J}(\tau-\mathcal T_0)^{\frac{1}{2}}\|\nabla T(\tau, \mathcal T_0)w(\mathcal T_0)\|_{3}\le C\|w(\mathcal T_0)\|_3
\end{align}
and from \eqref{dec-0-alpha1}--\eqref{dec-1-alpha1} that
\begin{align}\label{sob-wei-top}
&\sup_{\tau\in J}(\tau-\mathcal T_0)^{\frac{1}{2}}(1+\tau)^{-\alpha}\|\rho^{\alpha}T(\tau, \mathcal T_0)w(\mathcal T_0)\|_{\infty}+\sup_{\tau\in J}(\tau-\mathcal T_0)^{\frac{1}{2}}(1+\tau)^{-\alpha}\|\rho^{\alpha}\nabla T(\tau, \mathcal T_0)w(\mathcal T_0)\|_{3}\notag\\
\le &C\mathcal T_0^{-\alpha}\|\rho^{\alpha}w(\mathcal T_0)\|_{3}+C\|w(\mathcal T_0)\|_3,
\end{align}
where $J:=[\mathcal T_0, t]$.
From Remark \ref{rmk-const-j} together with \eqref{sol-rad}, we have
\begin{align}
&\left[\widetilde{H}[w, w]\right]_{\infty, 0, J}+\left[\nabla \widetilde{H}[w, w]\right]_{3, 0, J}\le 4C_0C_1^{\prime}\|\rho^{\alpha}w_0\|_3\left[\nabla w\right]_{3, 0, J},\label{sob-0-h}\\
&\left[\widetilde{G}[w]\right]_{\infty, 0, J}+\left[\nabla \widetilde{G}[w]\right]_{3, 0, J}\le C_2^{\prime}\llbracket v\rrbracket \Bigl(\left[\nabla w\right]_{3, 0, J}+\left[w\right]_{\infty, 0, J}\Bigr),\label{sob-0-g}
\end{align}
and 
\begin{align}
&\left[\widetilde{H}[w, w]\right]_{\infty, \alpha, J}+\left[\nabla \widetilde{H}[w, w]\right]_{3, \alpha, J}\le 8C_0C_1^{\prime\prime}\|\rho^{\alpha}w_0\|_3\left[\nabla w\right]_{3, 0, J},\label{sob-wei-h}\\
&\left[\widetilde{G}[w]\right]_{\infty, \alpha, J}+\left[\nabla \widetilde{G}[w]\right]_{3, \alpha, J}\le C_2^{\prime\prime}\llbracket v\rrbracket \Bigl(\left[\nabla w\right]_{3, 0, J}+\left[\nabla w\right]_{3, \alpha, J}+\left[w\right]_{\infty, 0, J}+\left[w\right]_{\infty, \alpha, J}\Bigr),\label{sob-wei-g}
\end{align}
where $[w]_{\infty, \sigma, J}$ and $[\nabla w]_{3, \sigma, J}$ are given by \eqref{sob-infty-t0}--\eqref{sob-nab-t0} respectively.
Collecting \eqref{sob-0-top}, \eqref{sob-0-h} and \eqref{sob-0-g} leads us to
\begin{align}\label{sob-0-preabs}
\left[w\right]_{\infty, 0, J}+\left[\nabla w\right]_{3, 0, J}\le C\|w(\mathcal T_0)\|_{3}+4C_0C_1^{\prime}\|\rho^{\alpha}w_0\|_3\left[\nabla w\right]_{3, 0, J}+C_2^{\prime}\llbracket v\rrbracket \Bigl(\left[\nabla w\right]_{3, 0, J}+\left[w\right]_{\infty, 0, J}\Bigr).
\end{align}
By taking into account \eqref{cond-small-per}, if the conditions
\begin{align}\label{cond-sob-1}
\|\rho^{\alpha}w_0\|_3\le\min\left\{\varepsilon_2^{\prime\prime}, \frac{1}{8C_0C_1^{\prime}}\right\}=:\varepsilon_2,\qquad \lambda^{-\frac{3}{4}}\|(\eta, \omega)\|_{W^{2, 2}(0, \mathcal T)}\le \min\left\{\varepsilon_1, \frac{1}{4C_2^{\prime}c_*}\right\}
\end{align}
are satisfied, it follows from \eqref{sob-0-preabs} that
\begin{align}\label{0-abs}
[w]_{\infty, 0, J}+[\nabla w]_{3, 0, J}\le 4C\|w(\mathcal T_0)\|_3,
\end{align}
where $\varepsilon_2^{\prime\prime}$ is given by \eqref{def-ep2pp}.
From this, we are led to
\begin{align*}
(t-\mathcal T_0)^{\frac{1}{2}}\|w(t)\|_{\infty}+(t-\mathcal T_0)^{\frac{1}{2}}\|\nabla w(\mathcal T_0)\|_{3}\le 4C\|w(\mathcal T_0)\|_3
\end{align*}
for every $t\ge\mathcal T_0$.
Since we consider the large time behavior, we may assume $t\ge 2\mathcal T_0$.
With the aid of the inequality $t^{\frac{1}{2}}\le \sqrt2(t-\mathcal T_0)^{\frac{1}{2}}$, we infer
\begin{align*}
t^{\frac{1}{2}}\|w(t)\|_{\infty}+t^{\frac{1}{2}}\|\nabla w(t)\|_{3}\le 4\sqrt2 C\|w(\mathcal T_0)\|_{3}
\end{align*}
for all $(\mathcal T_0, t)$ with $\mathcal T_0>0$ and $t\ge 2\mathcal T_0$.
Given $t>0$, we put $\mathcal T_0=t/2$ and take into account \eqref{sol-0-sob} with $q=3$ to find the former property of \eqref{sob-step2}.\par
It follows from \eqref{sob-wei-top} and \eqref{sob-wei-h}--\eqref{sob-wei-g} that
\begin{align}\label{sup-wei-sob-inf}
\begin{split}
[w]_{\infty, \alpha, J}+[\nabla w]_{3, \alpha, J}
&\le C\mathcal T_0^{-\alpha}\|\rho^{\alpha}w(\mathcal T_0)\|_{3}+C\|w(\mathcal T_0)\|_3+8C_0C_1^{\prime\prime}\|\rho^{\alpha}w_0\|_3\left[\nabla w\right]_{3, 0, J}\\
&+C_2^{\prime\prime}\llbracket v\rrbracket \Bigl(\left[\nabla w\right]_{3, 0, J}+\left[\nabla w\right]_{3, \alpha, J}+\left[w\right]_{\infty, 0, J}+\left[w\right]_{\infty, \alpha, J}\Bigr).
\end{split}
\end{align}
Suppose
\begin{align}\label{def-ep3}
\|\rho^{\alpha}w_0\|_3\le\varepsilon_2,\qquad \lambda^{-\frac{3}{4}}\|(\eta, \omega)\|_{W^{2, 2}(0, \mathcal T)}\le \min\left\{\varepsilon_1. \frac{1}{4C_2^{\prime}c_*}, \frac{1}{2C_2^{\prime\prime}c_*}\right\}=:\varepsilon_3,
\end{align}
where $\varepsilon_2$ is given by \eqref{cond-sob-1}.
Then, we see from \eqref{cond-small-per}, \eqref{0-abs} and \eqref{sup-wei-sob-inf} that
\begin{align*}
[\nabla w]_{3, \alpha, J}+[w]_{\infty, \alpha, J}
\le 2C\mathcal T_0^{-\alpha}\|\rho^{\alpha}w(\mathcal T_0)\|_3+\left(6C+8\frac{C_1^{\prime\prime}}{C_1^{\prime}}C\right)\|w(\mathcal T_0)\|_3.
\end{align*}
By the same reasoning as before, we find
\begin{align*}
t^{\frac{1}{2}-\alpha}\|\rho^{\alpha}w(t)\|_{\infty}+t^{\frac{1}{2}-\alpha}\|\rho^{\alpha}\nabla w(t)\|_3\le C\mathcal T_0^{-\alpha}\|\rho^{\alpha}w(\mathcal T_0)\|_3+C\|w(\mathcal T_0)\|_{3}
\end{align*}
for all $(\mathcal T_0, t)$ with $\mathcal T_0>0$ and $t\ge 2\mathcal T_0$.
Given $t>0$, taking $\mathcal T_0=t/2$ leads us to the latter property of \eqref{sob-step2} on account of \eqref{sol-wei-sob} with $q=3$.
\par
Finally, we verify \eqref{sol-0-sob}--\eqref{sol-wei-sob} for $r\in (3, \infty)$.
Let $t\ge 4$.
We use the equation
\begin{align*}
w(t)=T(t, t/2)w(t/2)-\int_{\frac{t}{2}}^{t}T(t, \tau)P(w\cdot\nabla w+v\cdot\nabla w+w\cdot\nabla v)(\tau)\,d\tau,
\end{align*}
and set
\begin{align*}
\langle w\rangle_{\infty, \sigma, t}:=\sup_{\frac{t}{2}<\tau<t}\tau^{\frac{1}{2}-\sigma}\|\rho^{\sigma}w(\tau)\|_{\infty},\qquad \langle \nabla w\rangle_{3, \sigma, t}:=\sup_{\frac{t}{2}<\tau<t}\tau^{\frac{1}{2}-\sigma}\|\rho^{\sigma}\nabla w(\tau)\|_{3}
\end{align*}
for $\sigma\in \{0, \alpha\}$.
Let $\beta:=\max\{0, 1/2-3/2r-\alpha/2\}$.
By \eqref{lqlr-hi-1} and \eqref{dec-1-alpha1}, we have
\begin{align}
&\left\|\nabla T\left(t, \frac{t}{2}\right)w\left(\frac{t}{2}\right)\right\|_{r}\le Ct^{-\frac{1}{2}}\left\|w\left(\frac{t}{2}\right)\right\|_3,\label{r-top-0}\\
&\left\|\rho^{\alpha}\nabla T\left(t, \frac{t}{2}\right)w\left(\frac{t}{2}\right)\right\|_r
\le Ct^{-1+\frac{3}{2r}+\beta}\left\|\rho^{\alpha}w\left(\frac{t}{2}\right)\right\|_3+Ct^{-1+\frac{3}{2r}+\beta+\alpha}\left\|w\left(\frac{t}{2}\right)\right\|_{3},\label{r-top-wei}
\end{align}
and
\begin{align}
&\int_{\frac{t}{2}}^{t}\|\nabla T(t, \tau)P(w\cdot\nabla w)(\tau)\|_r\,d\tau\le Ct^{-\frac{1}{2}}\|w\|_X\langle w\rangle_{\infty, 0, t},\label{r-nonlin-0}\\
&\int_{\frac{t}{2}}^{t}\|\rho^{\alpha}\nabla T(t, \tau)P(w\cdot\nabla w)(\tau)\|_r\,d\tau\le Ct^{-1+\frac{3}{2r}+\beta+\alpha}\|w\|_X\bigl(\langle w\rangle_{\infty, \alpha, t}+\langle w\rangle_{\infty, 0, t}\bigr).\label{r-nonlin-wei}
\end{align}
Following the same calculation as in Lemma \ref{stab-lem-duha}, we get
\begin{align*}
\int_{\frac{t}{2}}^{t-1}\|\nabla T(t, \tau)P(v\cdot\nabla w+w\cdot\nabla v)(\tau)\|_r\,d\tau
\le Ct^{-\frac{1}{2}}\llbracket v\rrbracket \bigl(\langle \nabla w\rangle_{3, 0, t}+\langle w\rangle_{\infty, 0, t}\bigr)
\end{align*}
and
\begin{align*}
&\int_{\frac{t}{2}}^{t-1}\|\rho^{\alpha}\nabla T(t, \tau)P(v\cdot\nabla w+w\cdot\nabla v)(\tau)\|_r\,d\tau\\
\le &Ct^{-\frac{1}{2}+\alpha}\int_{\frac{t}{2}}^{t-1}(t-\tau)^{-\frac{3}{2}-\frac{\delta_1}{2}+\frac{3}{2r}+\beta}\,d\tau\llbracket v\rrbracket\bigl(\langle \nabla w\rangle_{3, \alpha, t}+\langle \nabla w\rangle_{3, 0, t}\bigr)\\
+ &Ct^{-\frac{1}{2}+\alpha}\int_{\frac{t}{2}}^{t-1}(t-\tau)^{-\frac{3}{2}-\frac{\delta_2}{2}+\frac{3}{2r}+\beta}\,d\tau\llbracket v\rrbracket \bigl(\langle w\rangle_{\infty, \alpha, t}+\langle w\rangle_{\infty, 0, t}\bigr).
\end{align*}
On account of $-3/2-\delta_j/2+3/2r+\beta<-1$ for $j\in \{1, 2\}$, we get
\begin{align*}
&\int_{\frac{t}{2}}^{t-1}\|\nabla T(t, \tau)P(v\cdot\nabla w+w\cdot\nabla v)(\tau)\|_r\,d\tau\le Ct^{-\frac{1}{2}}\llbracket v\rrbracket\bigl(\langle \nabla w\rangle_{3, 0, t}+\langle w\rangle_{\infty, 0, t} \bigr),\\
&\int_{\frac{t}{2}}^{t-1}\|\rho^{\alpha}\nabla T(t, \tau)P(v\cdot\nabla w+w\cdot\nabla v)(\tau)\|_{r}\,d\tau\\
\le &Ct^{-\frac{1}{2}+\alpha}\llbracket v\rrbracket \bigl(\langle \nabla w\rangle_{3, \alpha, t}+\langle \nabla w\rangle_{3, 0, t}+\langle w\rangle_{\infty, \alpha, t}+\langle w\rangle_{\infty, 0, t}\bigr).
\end{align*}
As for the integral on $[t-1, t]$, estimate \eqref{est-uni-wei} leads us to
\begin{align*}
&\int_{t-1}^{t}\|\nabla T(t, \tau)P(v\cdot\nabla w+w\cdot\nabla v)(\tau)\|_{r}\,d\tau
\le Ct^{-\frac{1}{2}}\int_{t-1}^{t}(t-\tau)^{-1+\frac{3}{2r}}\,d\tau\,\llbracket v\rrbracket\bigl(\langle \nabla w\rangle_{3, 0, t}+\langle w\rangle_{\infty, 0, t}\bigr),\\
&\int_{t-1}^{t}\|\rho^{\alpha}\nabla T(t, \tau)P(v\cdot\nabla w+w\cdot\nabla v)(\tau)\|_{r}\,d\tau\\
\le &Ct^{-\frac{1}{2}+\alpha}\int_{t-1}^{t}(t-\tau)^{-1+\frac{3}{2r}}\,d\tau\,\llbracket v\rrbracket\bigl(\langle \nabla w\rangle_{3, \alpha, t}+\langle w\rangle_{\infty, \alpha, t}\bigr).
\end{align*}
The integrals in the right-hand side above are always convergent whenever $r<\infty$; thus, we have
\begin{align}
&\int_{\frac{t}{2}}^{t}\|\nabla T(t, \tau)P(v\cdot\nabla w+w\cdot\nabla v)(\tau)\|_{r}\,d\tau\le Ct^{-\frac{1}{2}}\llbracket v\rrbracket\bigl(\langle \nabla w\rangle_{3, 0, t}+\langle w\rangle_{\infty, 0, t}\bigr),\label{r-lin-0}\\
\begin{split}
&\int_{\frac{t}{2}}^{t}\|\rho^{\alpha}\nabla T(t, \tau)P(v\cdot\nabla w+w\cdot\nabla v)(\tau)\|_{r}\,d\tau\\
\le{} &Ct^{-\frac{1}{2}+\alpha}\llbracket v\rrbracket \bigl(\langle \nabla w\rangle_{3, \alpha, t}+\langle \nabla w\rangle_{3, 0, t}+\langle w\rangle_{\infty, \alpha, t}+\langle w\rangle_{\infty, 0, t}\bigr).\label{r-lin-wei}
\end{split}
\end{align}
Collecting \eqref{r-top-0}--\eqref{r-lin-wei} leads us to
\begin{align}\label{nab-r-0}
&t^{\frac{1}{2}}\|w(t)\|_r\le C\left\|w\left(\frac{t}{2}\right)\right\|_3+C\|w\|_X\langle w\rangle_{\infty, 0, t}+C\llbracket v\rrbracket\bigl(\langle \nabla w\rangle_{3, 0, t}+\langle w\rangle_{\infty, 0, t}\bigr)
\end{align}
and
\begin{align}\label{nab-r-wei}
\begin{split}
t^{\frac{1}{2}-\alpha}\|\rho^{\alpha}\nabla w(t)\|_r
\le  &Ct^{-\frac{1}{2}+\frac{3}{2r}+\beta}\left(t^{-\alpha}\left\|\rho^{\alpha}w\left(\frac{t}{2}\right)\right\|_3+\left\|w\left(\frac{t}{2}\right)\right\|_{3}\right)\\
+&Ct^{-\frac{1}{2}+\frac{3}{2r}+\beta}\|w\|_X\bigl(\langle w\rangle_{\infty, \alpha, t}+\langle w\rangle_{\infty, 0, t}\bigr)\\
+&C\llbracket v\rrbracket \bigl(\langle \nabla w\rangle_{3, \alpha, t}+\langle \nabla w\rangle_{3, 0, t}+\langle w\rangle_{\infty, \alpha, t}+\langle w\rangle_{\infty, 0, t}\bigr).
\end{split}
\end{align}
Suppose \eqref{def-ep3}.
From \eqref{sol-0-sob}--\eqref{sol-wei-sob} for $q=3$ together with \eqref{sob-step2} and $-1/2+3/2r+\beta\le0$, all the terms in the right-hand side of \eqref{nab-r-0}--\eqref{nab-r-wei} converge to zero as $t\to\infty$. 
The proof is complete.
\hfill
$\Box$

\bigskip
\noindent
{\it Acknowledgements}.
I would like to thank Professor Toshiaki Hishida for valuable discussions to improve this paper.
I am also grateful to Professor Tomoki Takahashi for sending me his paper \cite{T26} before publication.
This work is financially supported by the THERS Make New Standards Program for the Next Generation Researchers
through JPMJSP2125 from JST.

\end{document}